\documentclass[twoside,11pt]{article}

\usepackage{jmlr2e}
\usepackage[T1]{fontenc}
\usepackage{ae}
\usepackage{aecompl}
\usepackage[ansinew]{inputenc}

\expandafter\let\csname equation*\endcsname\relax % Avoids a stupid clash between iopart and amsmath
\expandafter\let\csname endequation*\endcsname\relax % Avoids a stupid clash between iopart and amsmath
\usepackage{amsmath}  % Prettier equation layouts
\usepackage{amssymb}
\usepackage{stmaryrd} 
\usepackage{rotating}
\usepackage{graphicx}
\usepackage{subfig}
\usepackage{calc}
\usepackage{dsfont}
\usepackage{latexsym}
\usepackage{mathtools}
\usepackage{algorithmic}
\usepackage{epstopdf}
\DeclareGraphicsExtensions{.pdf,.eps,.png,.jpg,.mps}
\usepackage{epsfig}
\usepackage{url}
\newcommand{\mc}{\mathcal}            % Graph font
\newcommand{\G}{\mathcal{G}}   % Digraph symbol
\newcommand{\W}{\mathcal{W}}  % Walk set
\newcommand{\V}{\mathcal{V}}
\newcommand{\E}{\mathcal{E}}

\newcommand{\C}{\mathcal{C}}

\newcommand{\N}{\mathcal{N}}
\newcommand{\X}{\mathbf{X}}

\usepackage{color}

\newcommand\ci{\perp\!\!\!\perp}
\newcommand{\PathOnePentagon}{{\,\atop \mathord{\includegraphics[height=7ex]{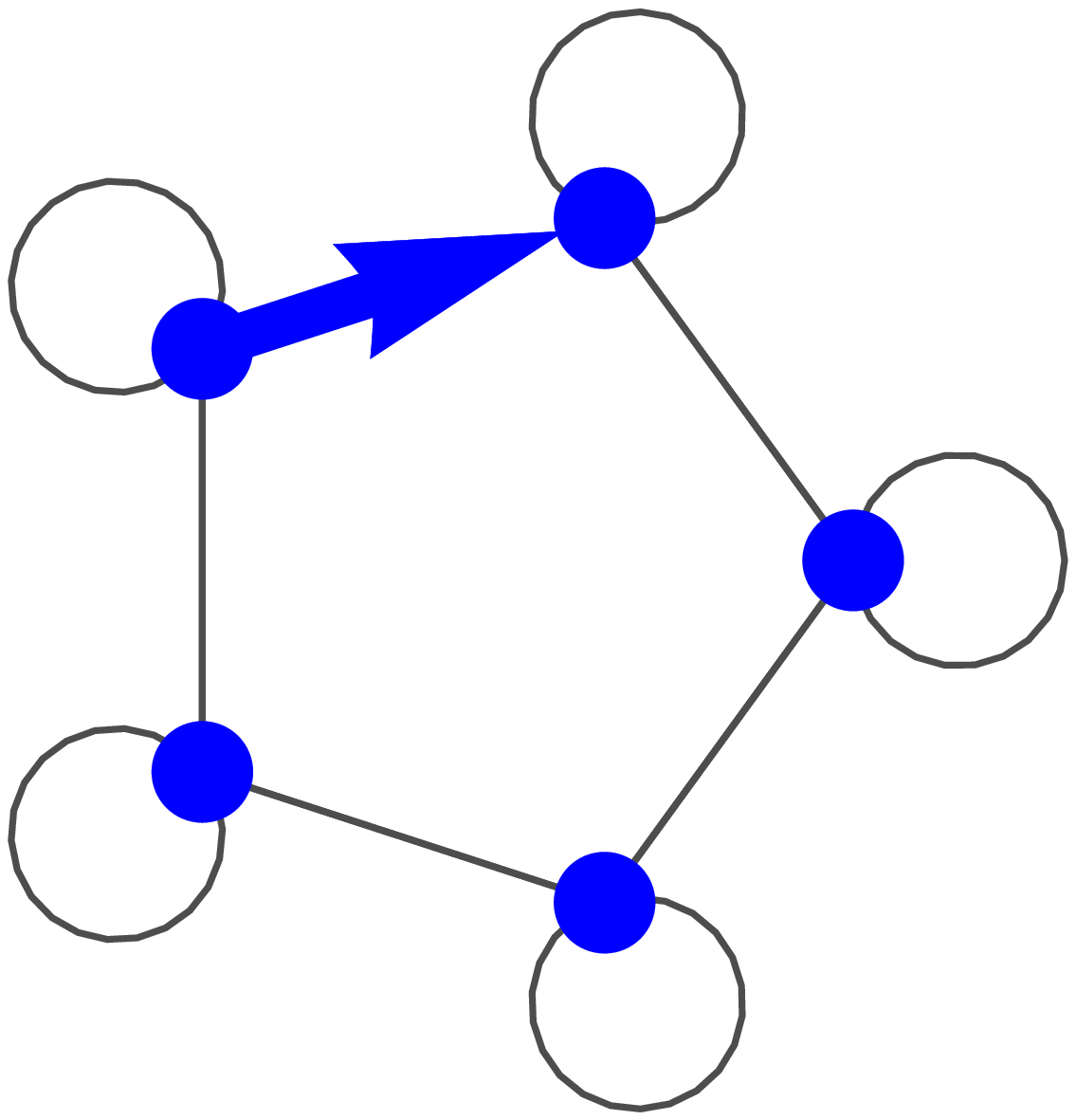}}}}
\newcommand{\PathTwoPentagon}{{\,\atop \mathord{\includegraphics[height=7ex]{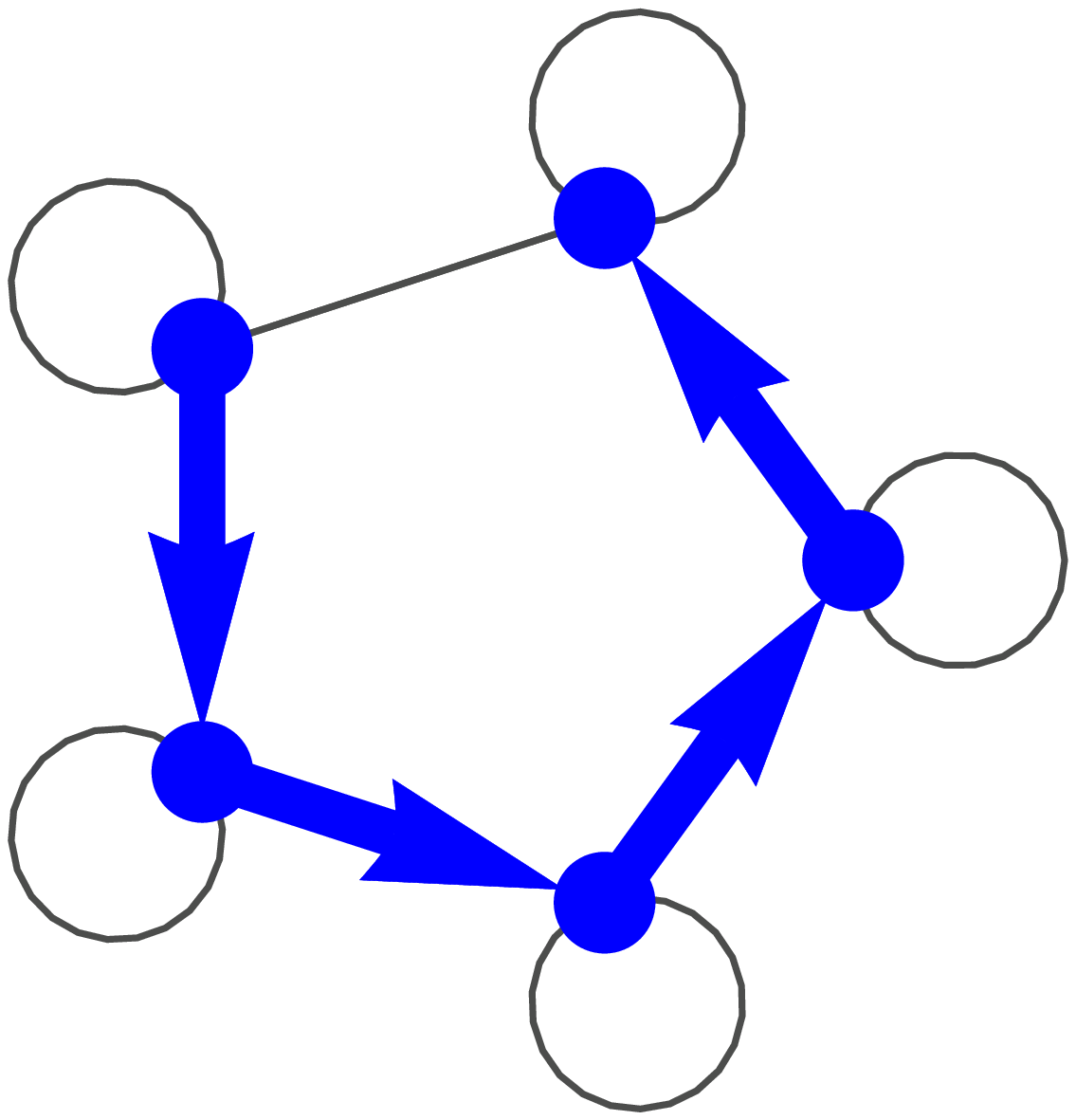}}}}
\newcommand{\SelfLoopThree}{{\,\atop \mathord{\includegraphics[height=7ex]{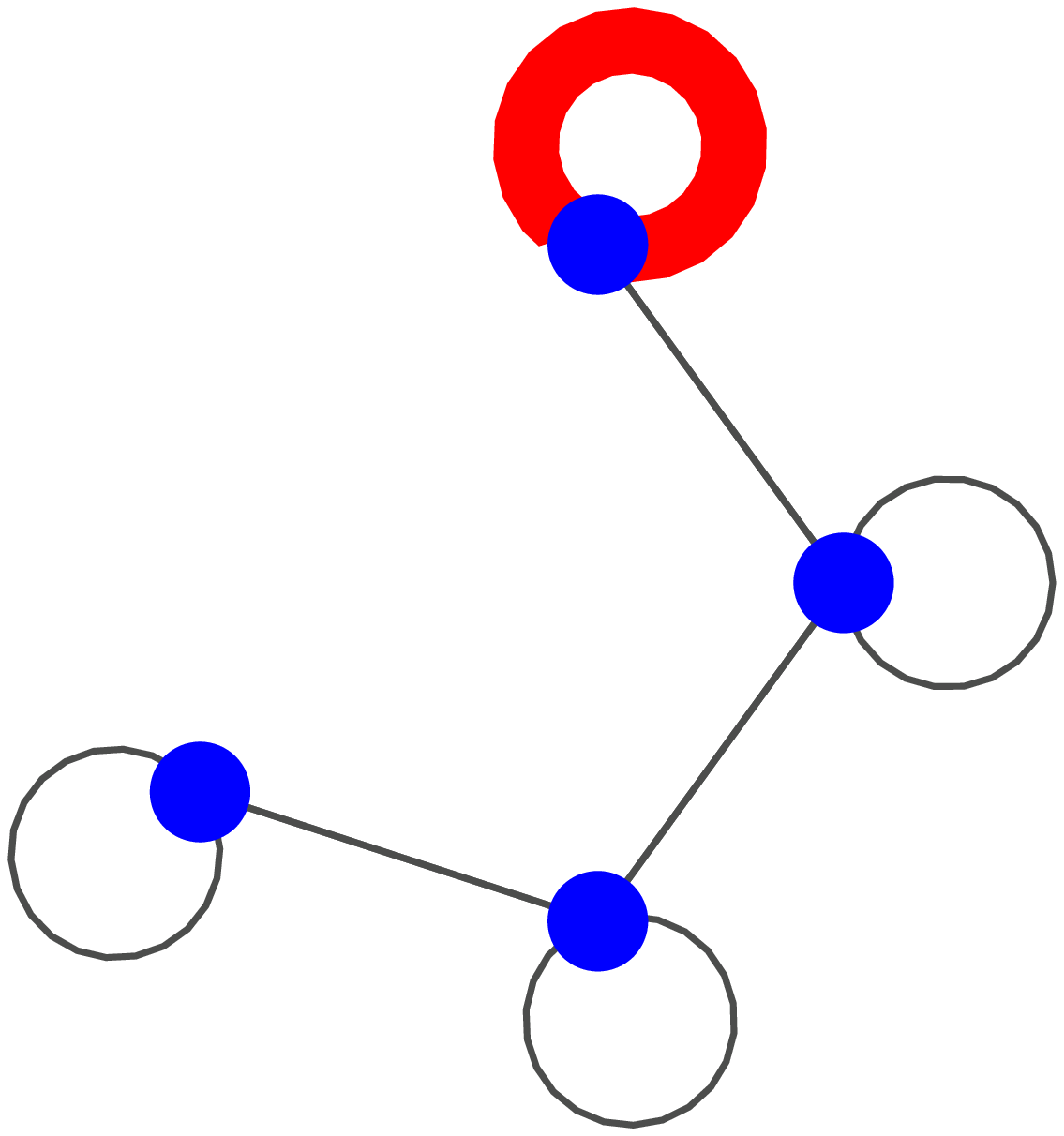}}}}
\newcommand{\BacktrackThree}{{\,\atop \mathord{\includegraphics[height=7ex]{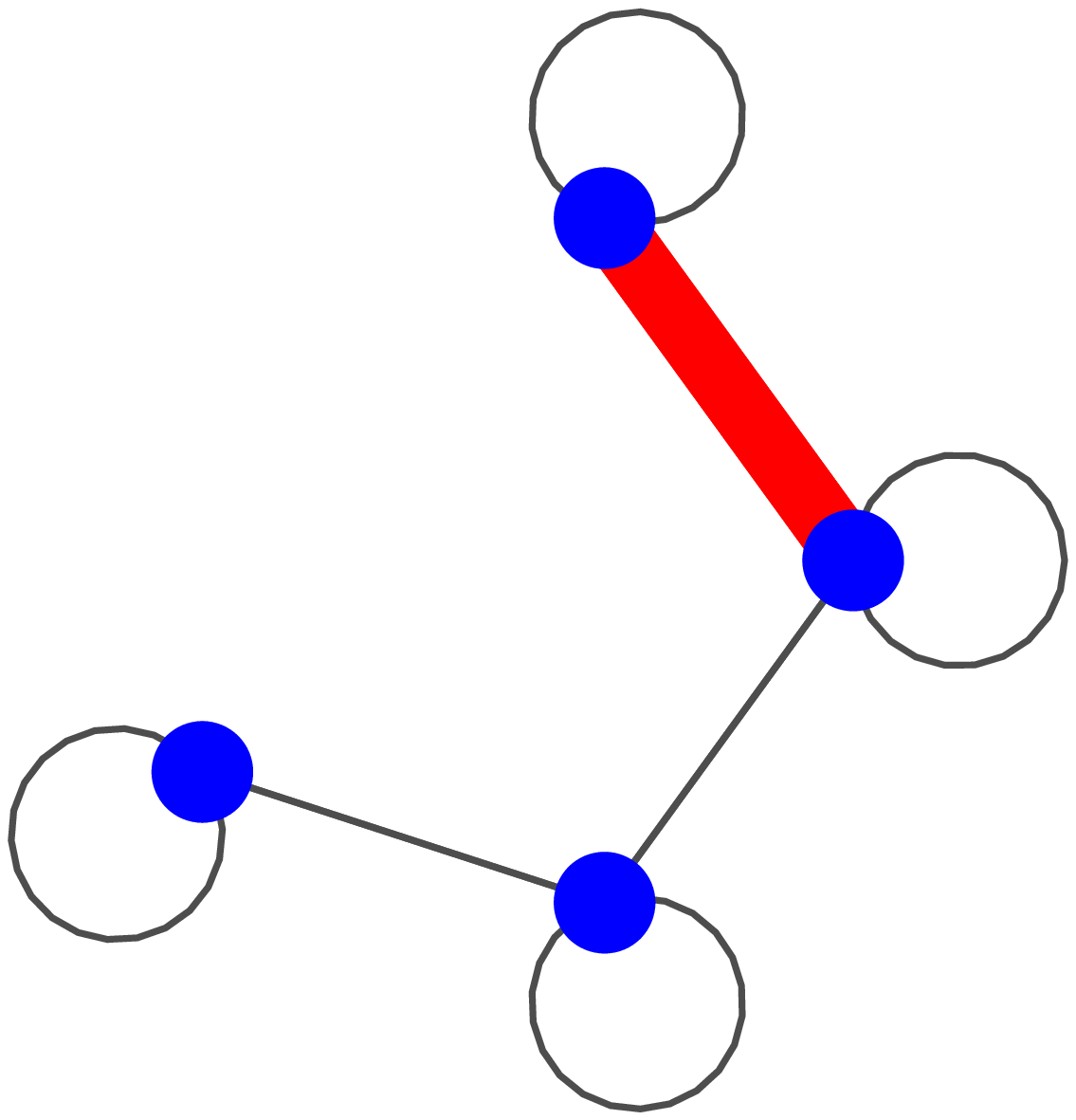}}}}
\newcommand{\SelfLoopOne}{{\,\atop \mathord{\includegraphics[height=7ex]{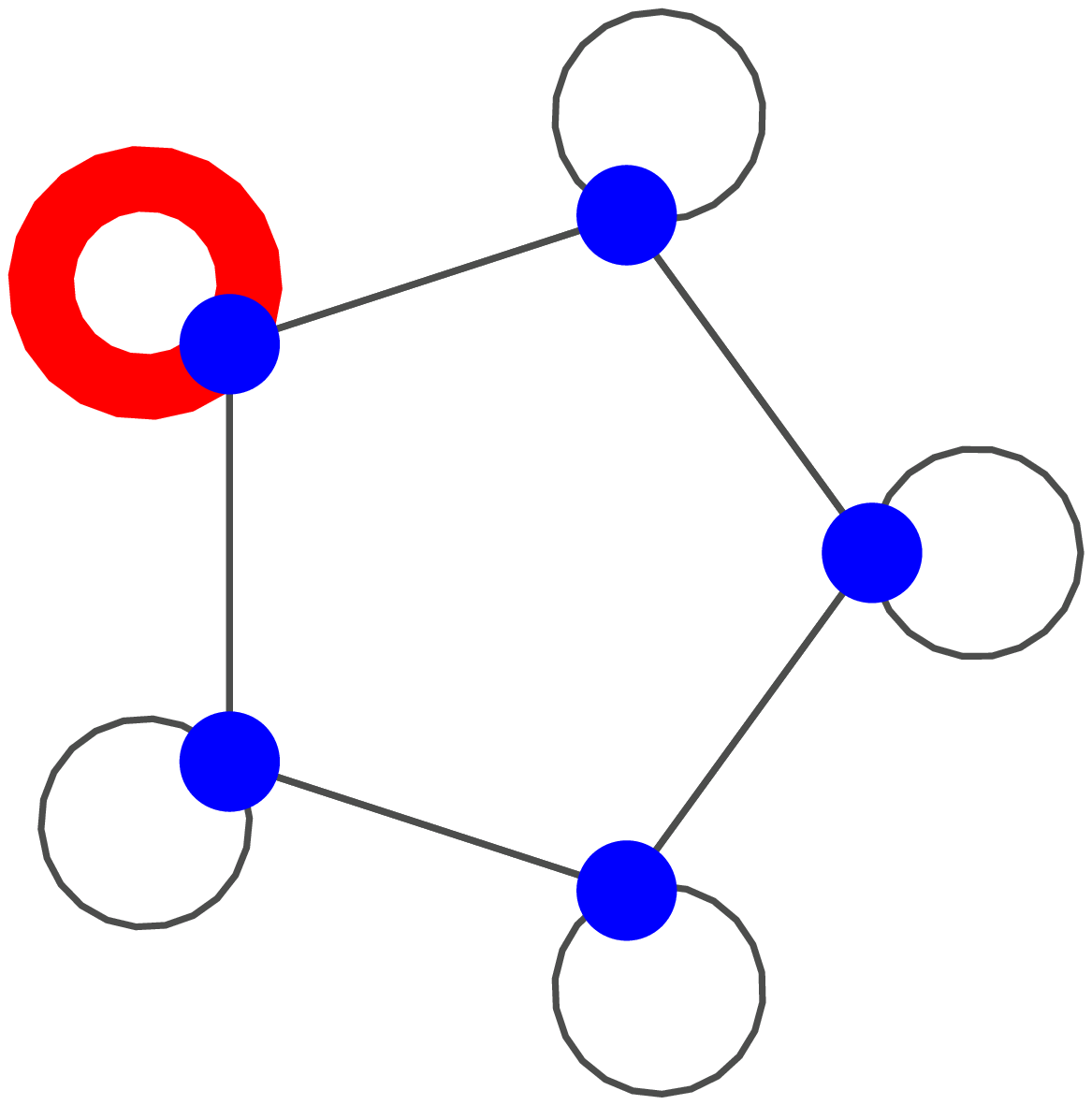}}}}
\newcommand{\BacktrackOne}{{\,\atop \mathord{\includegraphics[height=7ex]{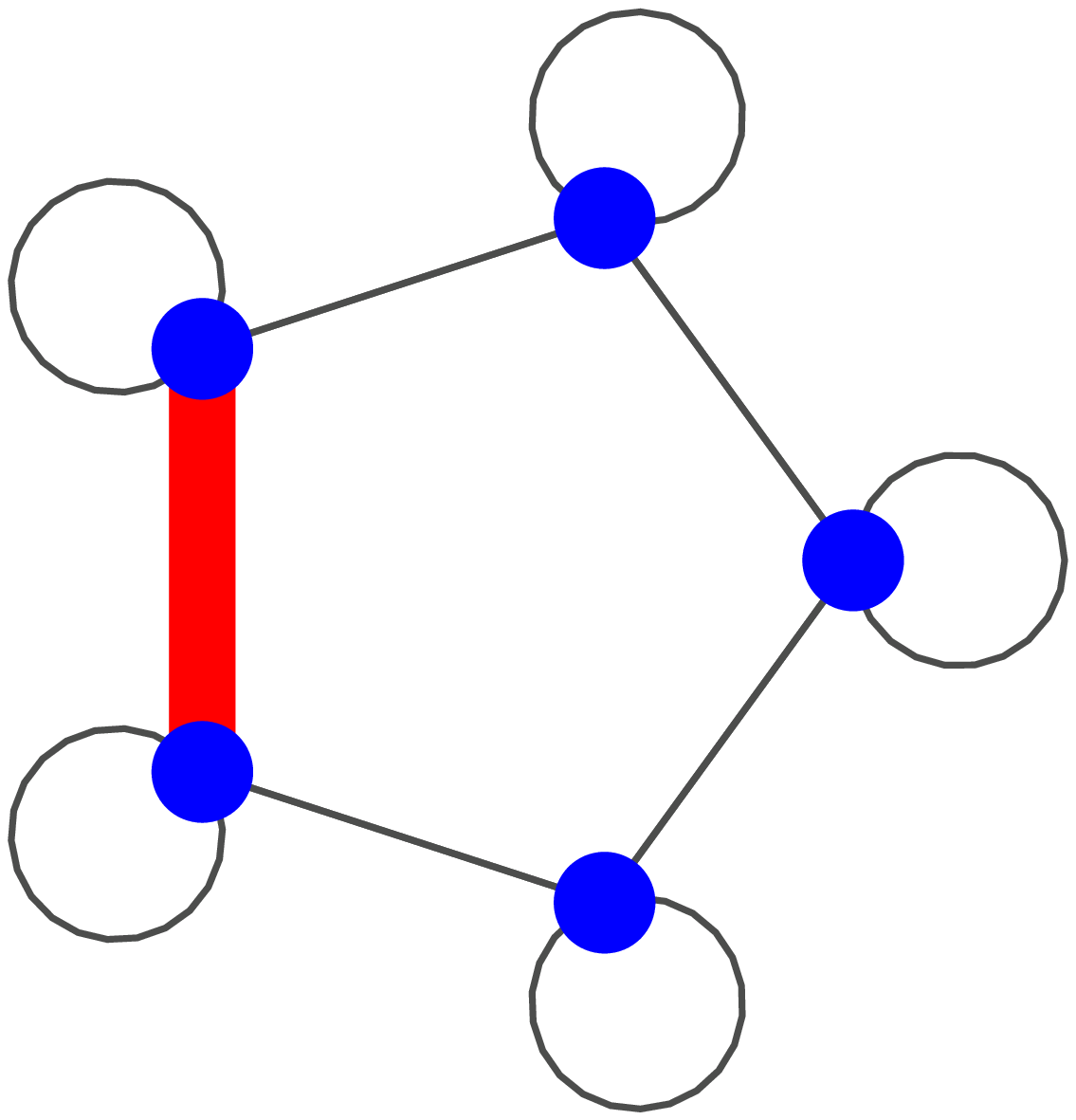}}}}
\newcommand{\BacktrackTwo}{{\,\atop \mathord{\includegraphics[height=7ex]{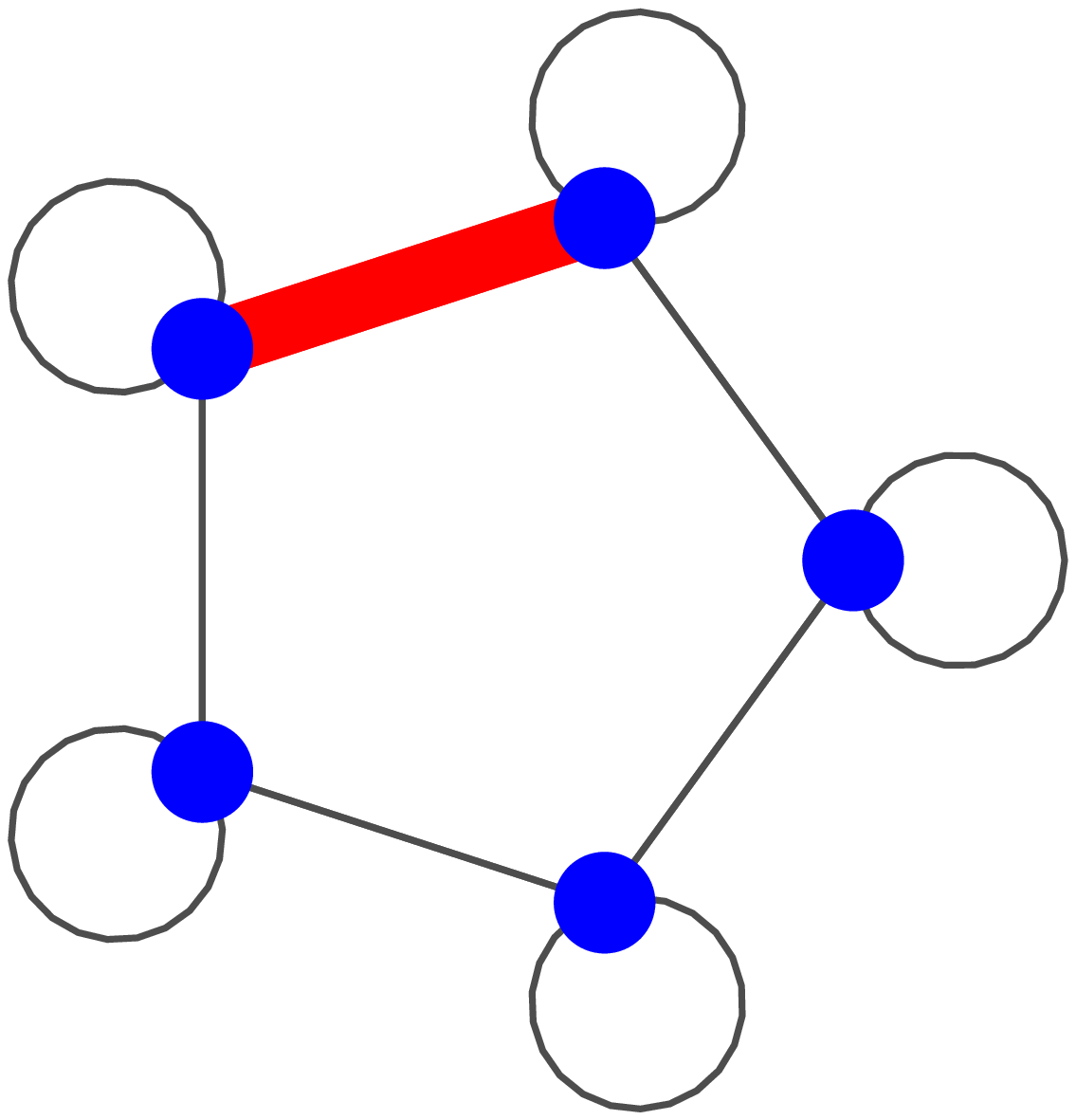}}}}
\newcommand{\PentagonOne}{{\,\atop \mathord{\includegraphics[height=7ex]{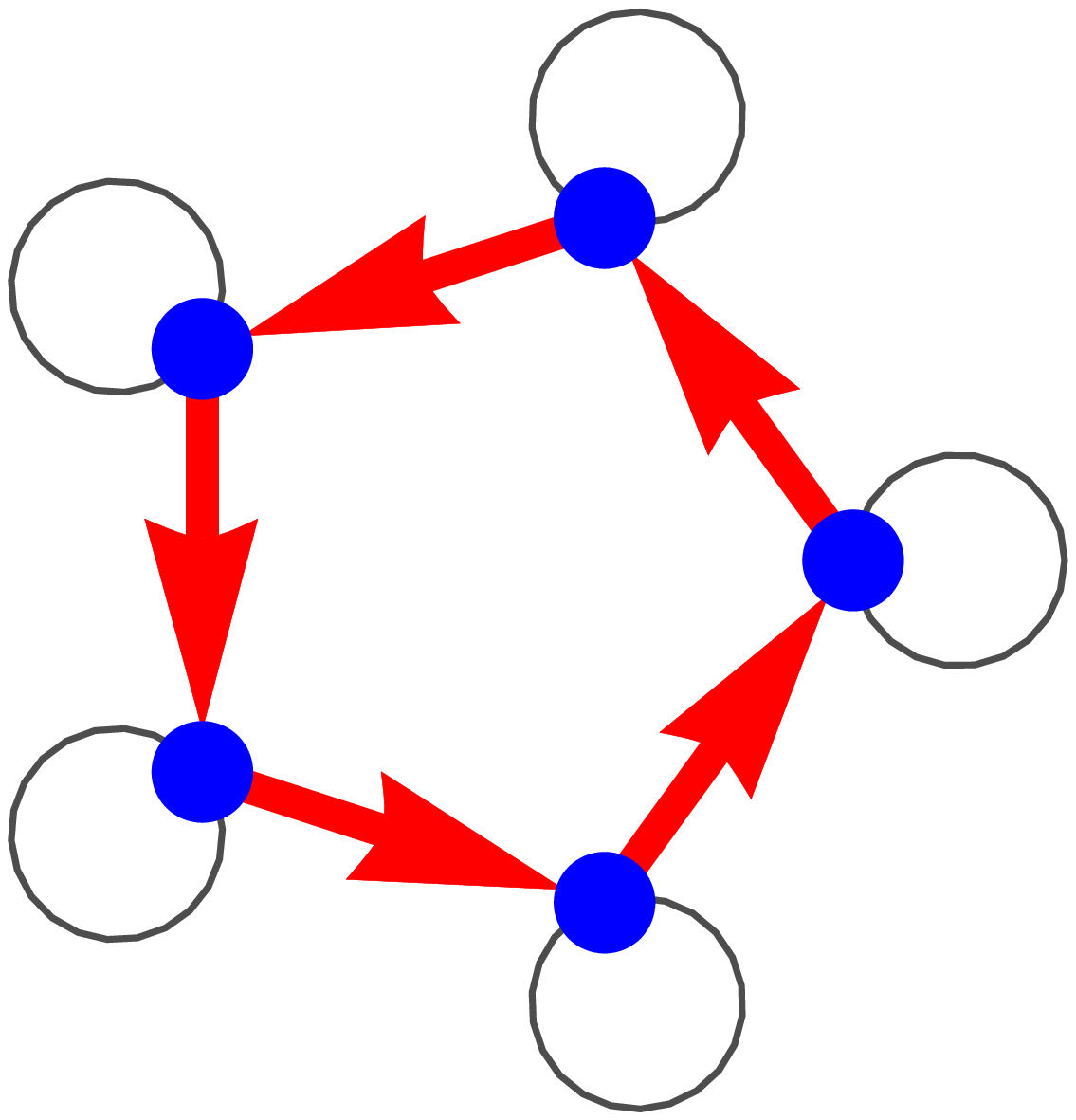}}}}
\newcommand{\PentagonTwo}{{\,\atop \mathord{\includegraphics[height=7ex]{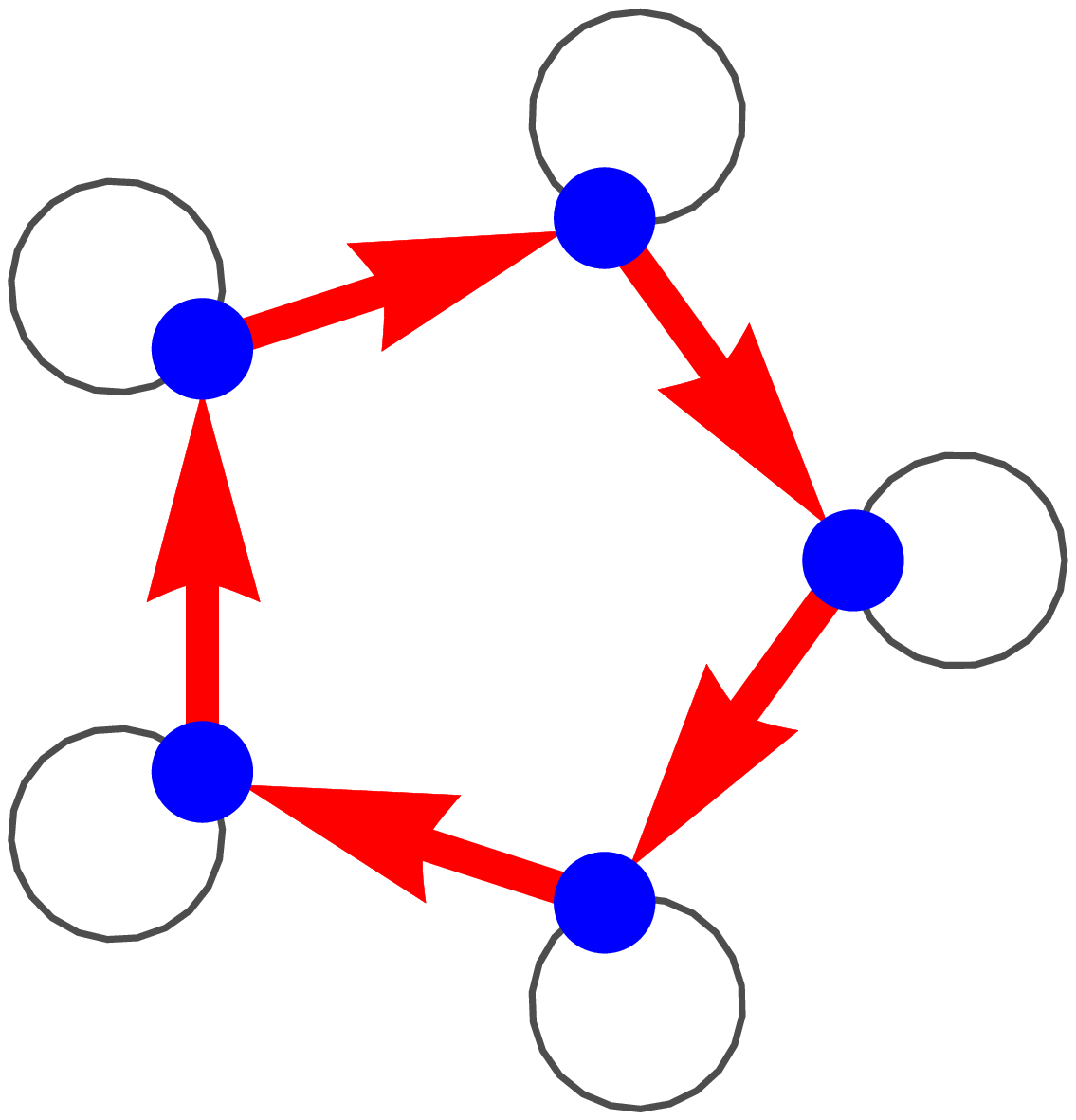}}}}

\newcommand{\LoopTri}{{\,\atop \mathord{\includegraphics[height=7ex]{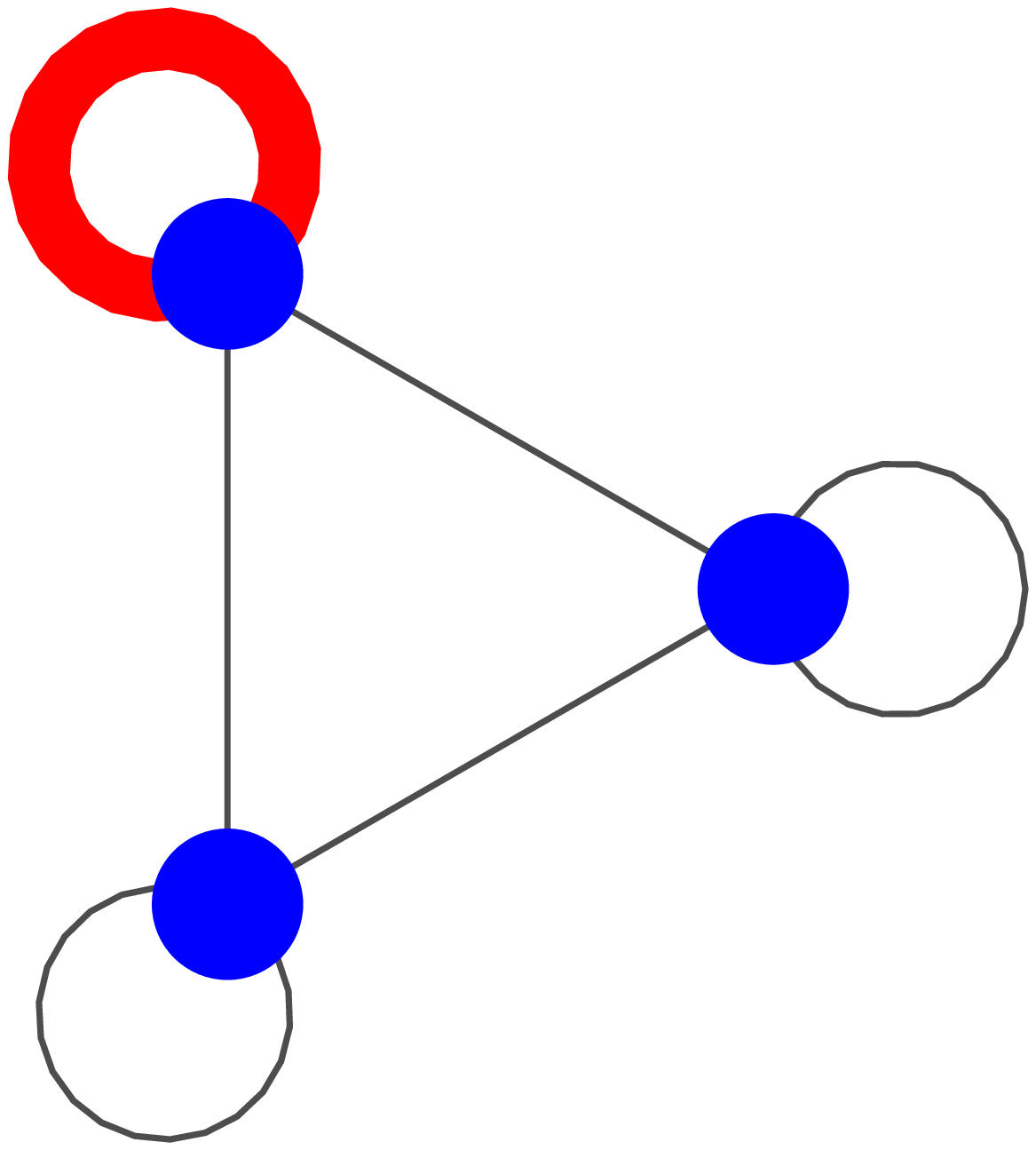}}}}
\newcommand{\BackTriOne}{{\,\atop \mathord{\includegraphics[height=7ex]{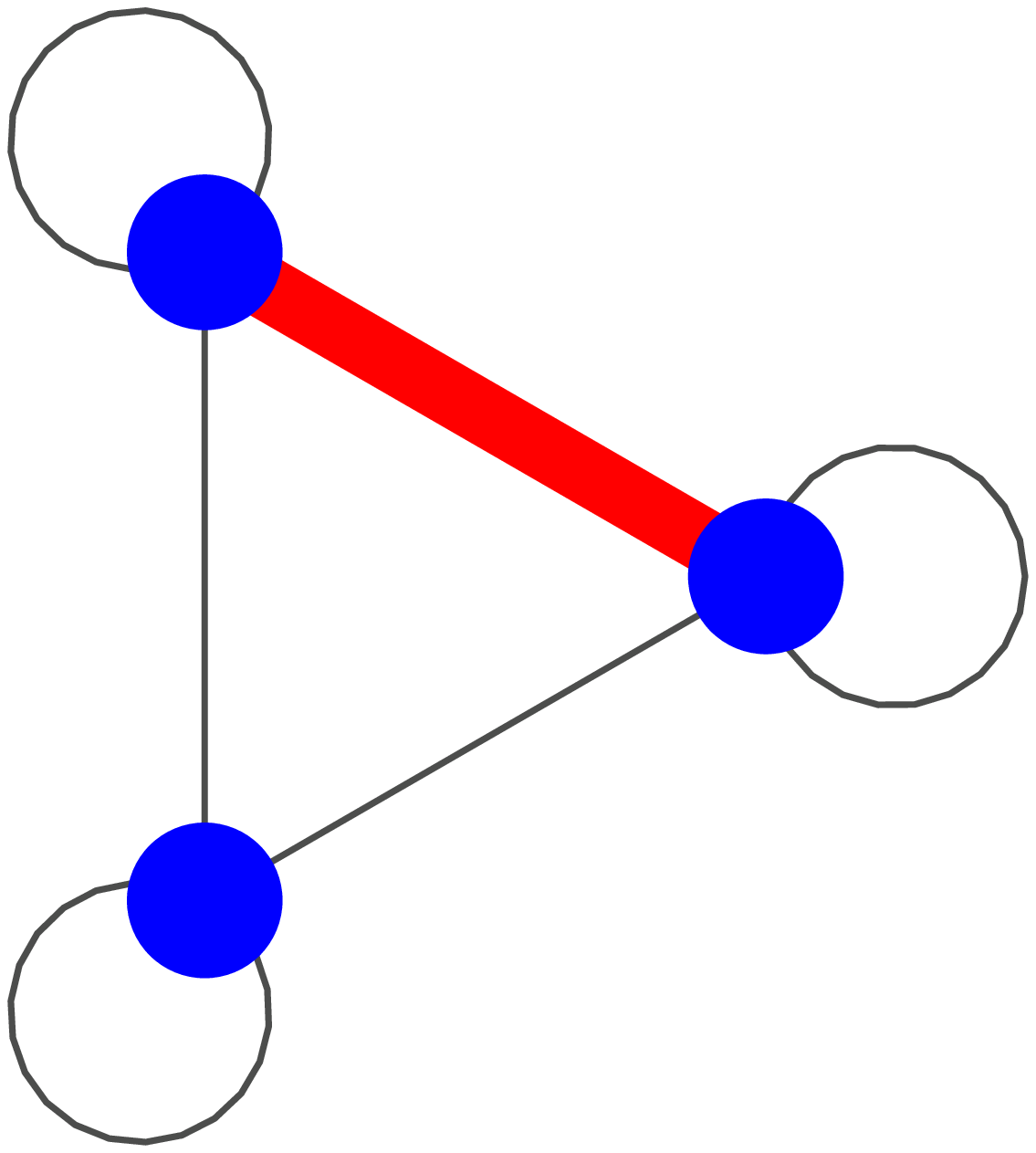}}}}
\newcommand{\BackTriTwo}{{\,\atop \mathord{\includegraphics[height=7ex]{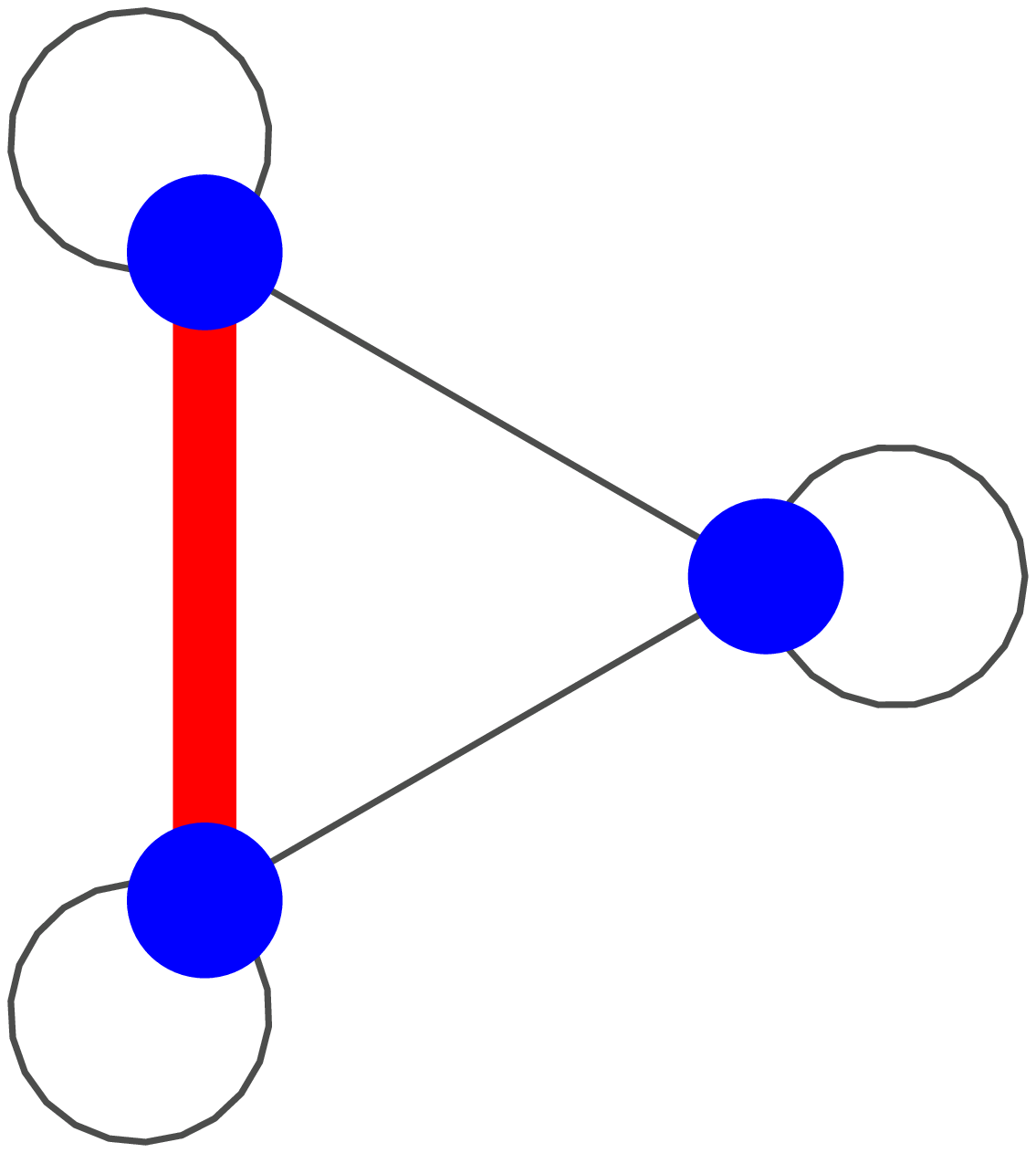}}}}
\newcommand{\TriOne}{{\,\atop \mathord{\includegraphics[height=7ex]{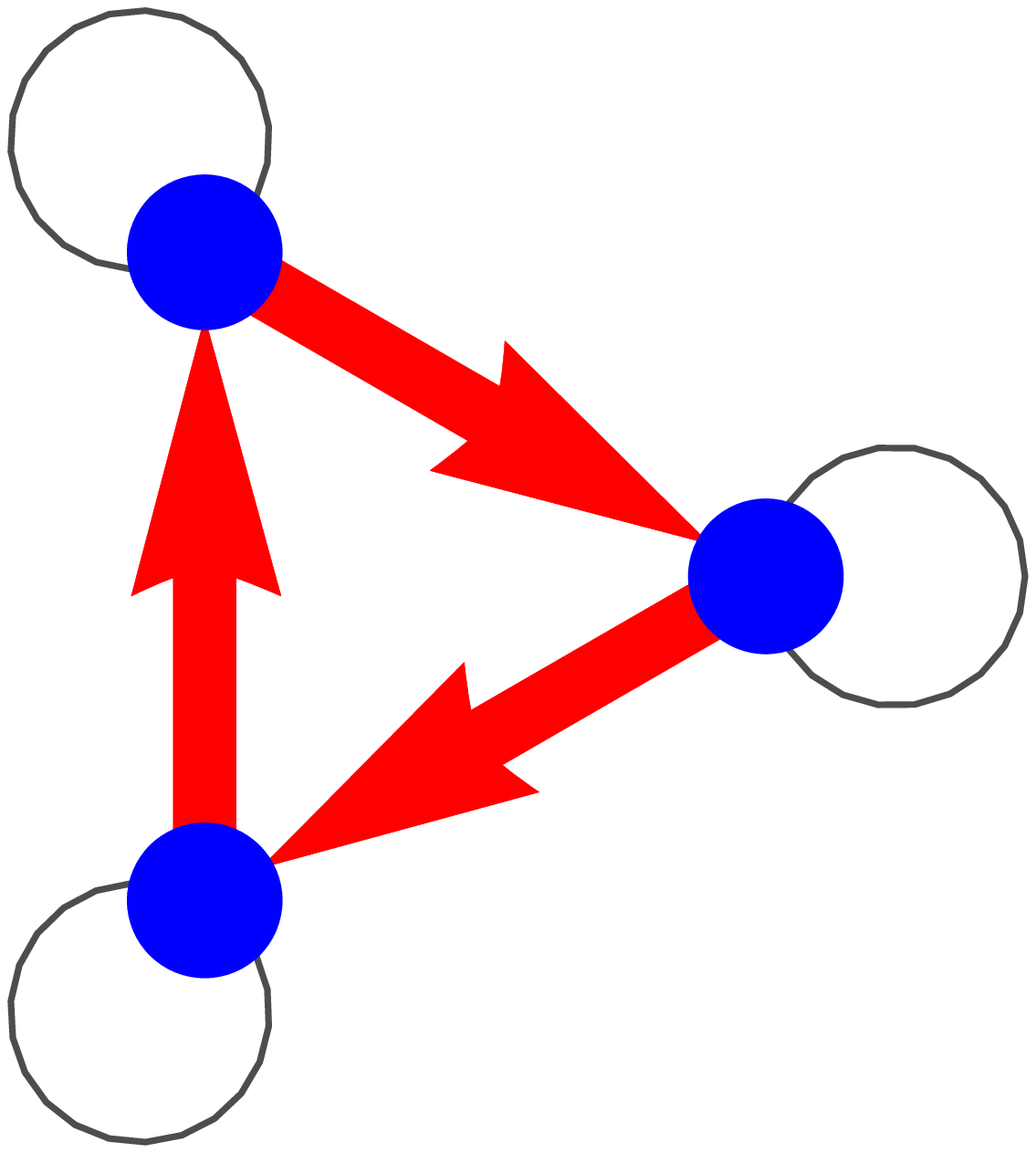}}}}
\newcommand{\TriTwo}{{\,\atop \mathord{\includegraphics[height=7ex]{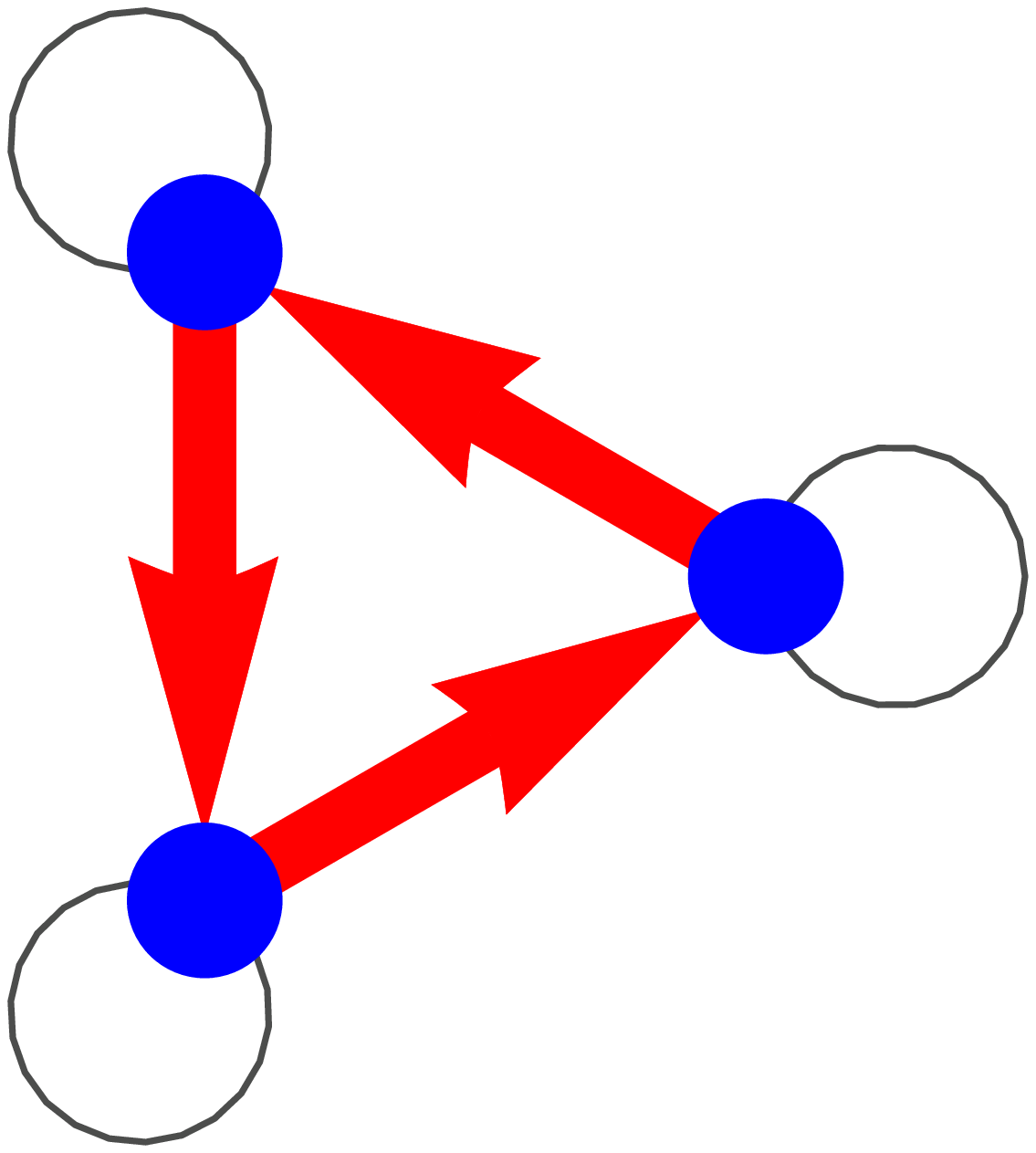}}}}

\newcommand{\LoopTriNoAlpha}{{\,\atop \mathord{\includegraphics[height=7ex]{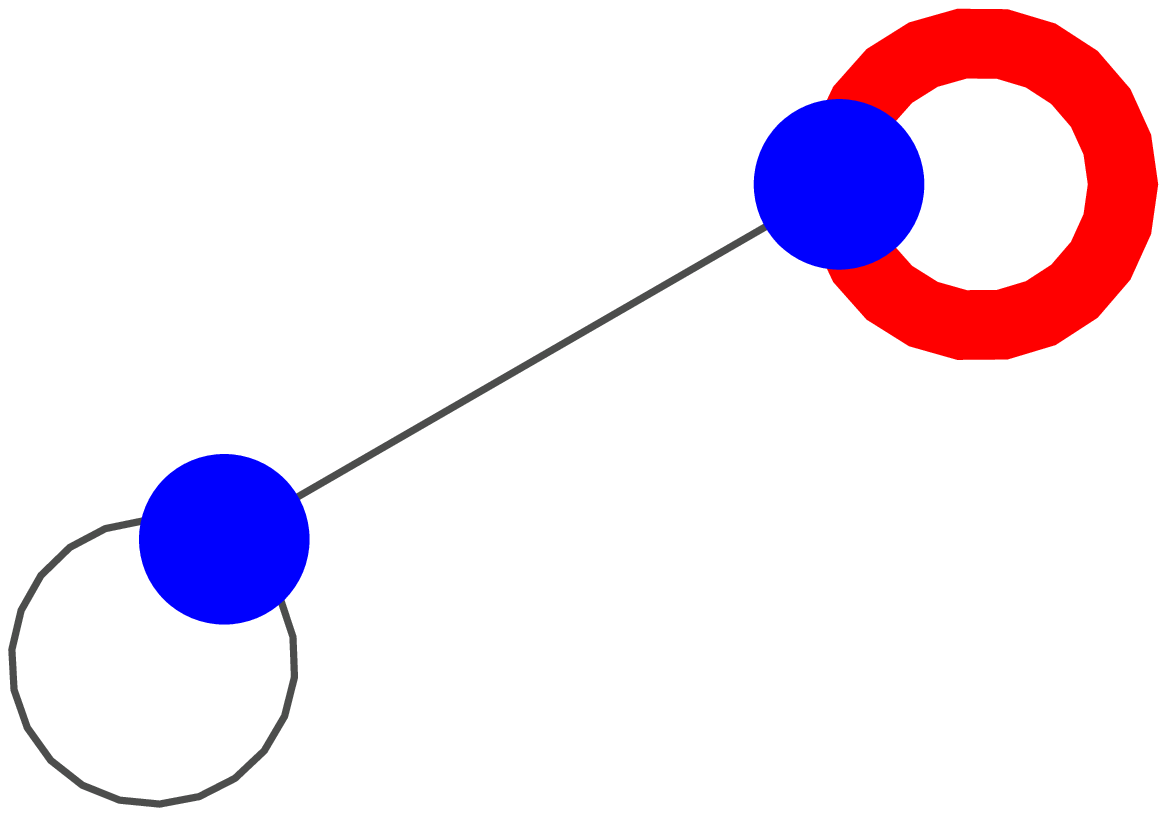}}}}
\newcommand{\BackTriNoAlpha}{{\,\atop \mathord{\includegraphics[height=7ex]{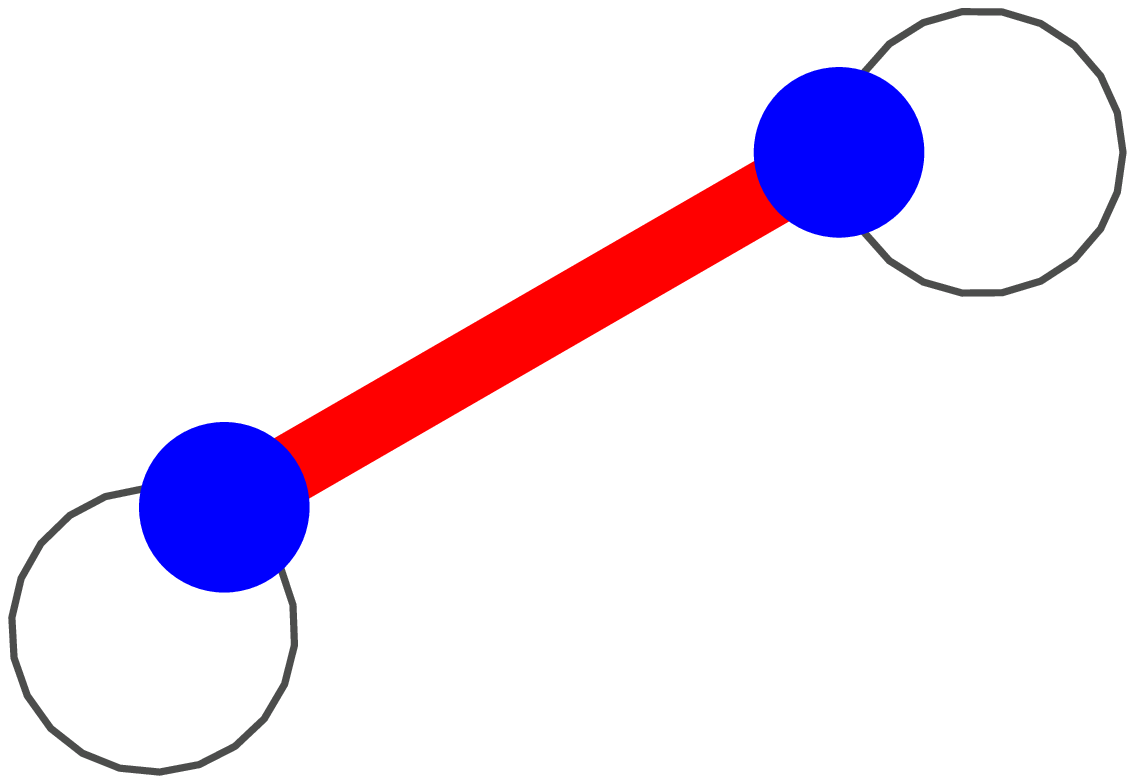}}}}
\newcommand{\LoopLoop}{{\,\atop \mathord{\includegraphics[height=7ex]{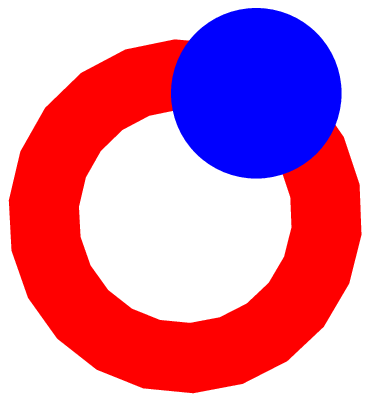}}}}

\newcommand{\PathOneTri}{{\,\atop \mathord{\includegraphics[height=7ex]{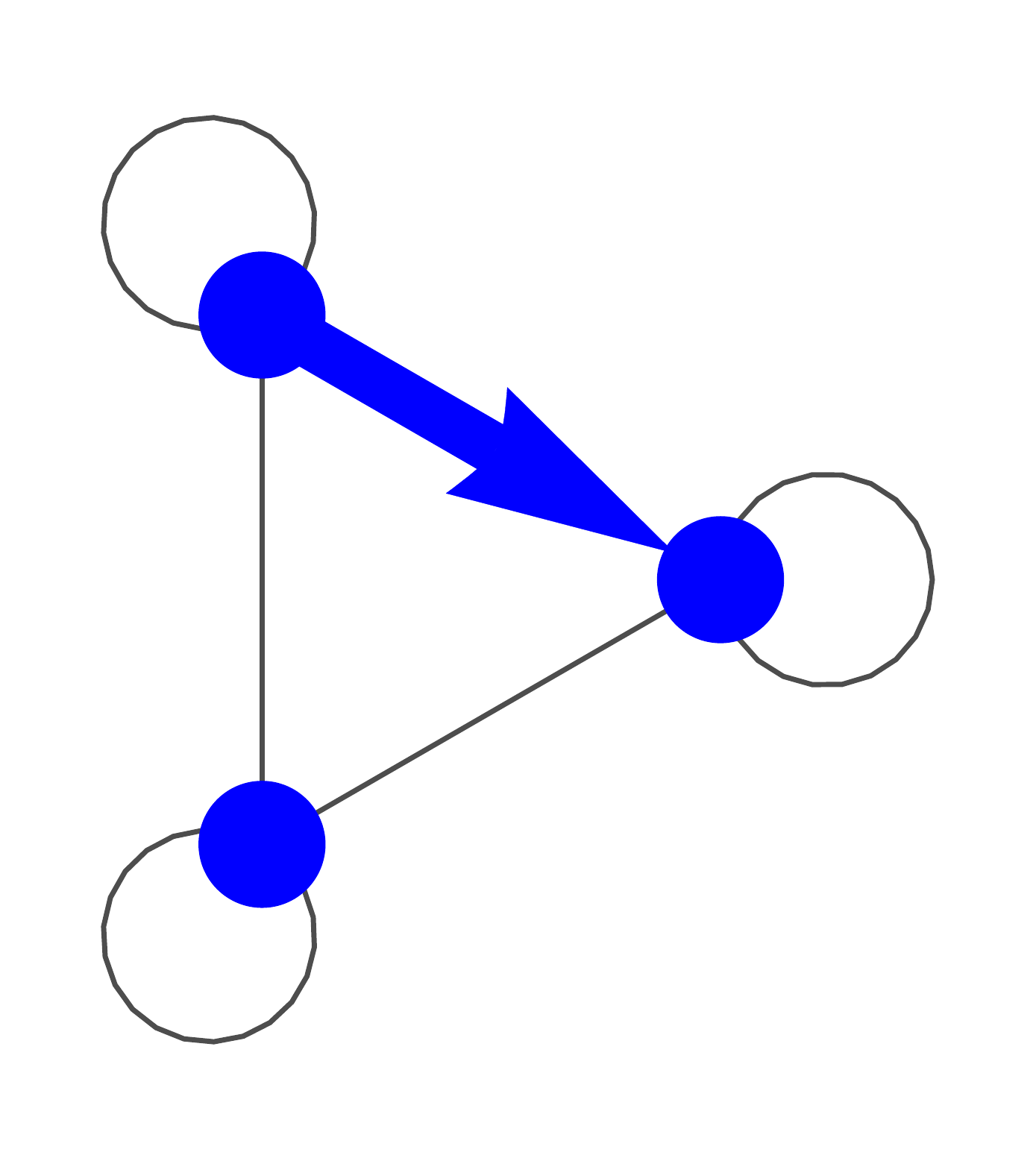}}}}
\newcommand{\PathTwoTri}{{\,\atop \mathord{\includegraphics[height=7ex]{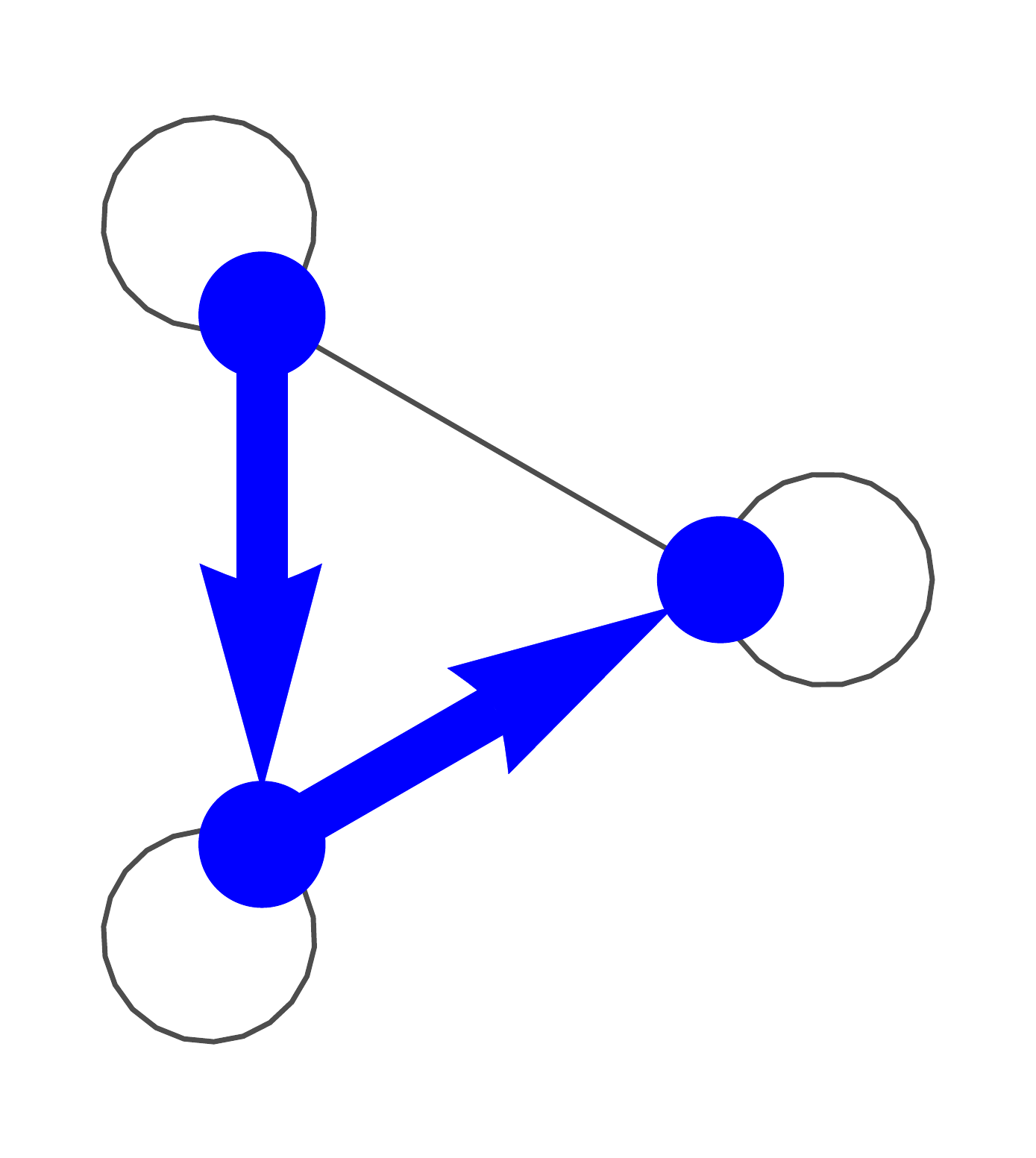}}}}

\jmlrheading{}{}{}{10/14}{}{P.-L. Giscard, Z. Choo, S. J. Thwaite and D. Jaksch}

\ShortHeadings{Exact Inference on Gaussian Graphical Models Using Path-Sums}{Giscard, Choo, Thwaite, and Jaksch}
\firstpageno{1}

\begin{document}

\title{Exact Inference on Gaussian Graphical Models of Arbitrary Topology using Path-Sums}
\author{\name P.-L.~Giscard \email p.giscard1@physics.ox.ac.uk \\
       \addr Department of Physics\\
       University of Oxford\\
       Clarendon Laboratory, Oxford OX1 3PU, United Kingdom
       \AND
       \name Z.~Choo  \email zheng.choo@stats.ox.ac.uk\\
       \addr Department of Statistics\\
       University of Oxford\\
       1 South Parks Road, Oxford OX1 3TG, UK
       \AND
       \name S.~J.~Thwaite \email simon.thwaite@physik.uni-muenchen.de\\
       \addr Department of Physics and Arnold Sommerfeld Center for Theoretical Physics,\\
Ludwig-Maximilians-Universit\"{a}t M\"{u}nchen\\ 
Theresienstra{\ss}e~37, 80333 Munich, Germany
       \AND
       \name D.~Jaksch \email d.jaksch1@physics.ox.ac.uk\\
       \addr Department of Physics\\
       University of Oxford\\
       Clarendon Laboratory, Oxford OX1 3PU, United Kingdom.\\
       \addr Centre for Quantum Technologies\\ 
       National University of Singapore\\ 3 Science Drive 2, Singapore 117543
       }

\editor{}

\date{\today}

\maketitle
\begin{abstract}
\hspace{-1.15mm}We present the path-sum formulation for exact statistical inference of marginals on Gaussian graphical models of arbitrary topology. The path-sum formulation gives the covariance between each pair of variables as a branched continued fraction of finite depth and breadth. Our method originates from the closed-form resummation of infinite families of terms of the walk-sum representation of the covariance matrix. We prove that the path-sum formulation always exists for models whose covariance matrix is positive definite: i.e.~it is valid for both walk-summable and non-walk-summable graphical models of arbitrary topology. We show that for graphical models on trees the path-sum formulation is equivalent to Gaussian belief propagation. We also recover, as a corollary, an existing result that uses determinants to calculate the covariance matrix. We show that the path-sum formulation formulation is valid for arbitrary partitions of the inverse covariance matrix. We give detailed examples demonstrating our results. 
\end{abstract}


\section{Introduction}
\label{sec:intro}
%
%
%\subsection{Background on Gaussian graphical models}
A Gaussian Markov random field (GMRF) is a random vector that follows a multivariate normal (or Gaussian) distribution and satisfies \emph{conditional independence} assumptions, hence the \emph{Markov} property. If $X_1,X_2,X_3$ are random variables with a joint probability density function (or joint probability mass function in a discrete case), we say that $X_1$ is \emph{conditionally independent} of $X_2$ given $X_3$, denoted $X_1\ci X_2|X_3$, if \citep{Lau:1996}
\begin{align*}
f(x_1,x_2|x_3)=f(x_1|x_3)f(x_2|x_3).
\end{align*}
Here we use $f$ as a generic symbol for the probability density function of the random variables corresponding to its arguments. GMRFs have a simple interpretation and find their applications, for example, in image analysis, spatial statistics, structural time series analysis and analysis of longitudinal and survival data \citep{RH:2005}. 

Consider a random vector $\X=(X_1,X_2,\dotsc,X_n)$ following a multivariate normal distribution with mean $\mu$ and covariance matrix $\Sigma$, denoted $\X\sim\N(\mu,\Sigma)$. The probability density function of $\X$ is given as
\begin{align*}
f(x)=\frac{1}{\sqrt{(2\pi)^n\det(\Sigma)}}\exp\left[-\frac{1}{2}(x-\mu)^T\Sigma^{-1}(x-\mu)\right].
\end{align*}    
Here $\Sigma$ is a symmetric and positive definite matrix. We write $A\succ0$ to denote that the matrix $A$ is positive definite. Alternatively the probability density function of $\X$ can be expressed in a \emph{canonical form}
\begin{align}
f(x)&=g(x)\exp\left[-\frac{1}{2}x^TJx+h^Tx-k(\mu,\Sigma)\right]\nonumber\\
&\propto\exp\left[-\frac{1}{2}x^TJx+h^Tx\right]\label{eq:cf},
\end{align}
where $J=\Sigma^{-1},h=J\mu,g(x)=(2\pi)^{-\frac{n}{2}}$ and $k(\mu,\Sigma)=-\frac{1}{2}\mu^TJ\mu-\frac{1}{2}\ln(\det(J))$. We call $J$ the \emph{information matrix} (or precision matrix) and $h$ the \emph{potential vector}. 

One advantage of using the form of parametrization in (\ref{eq:cf}) is that $J$ admits a graphical model in the following sense. Let $\G=(\V,\E)$ be an undirected graph with the vertex set $\V$ and the set of edges $\E$. Let $\X_{\setminus{ij}}$ denote the set of variables with $X_i$ and $X_j$ removed from $\X$. If $J=(J_{ij})_{i,j\in\V}$ is positive definite, then for $i,j\in\V$, where $i\neq j$, we have (Proposition 5.2, \citealp{Lau:1996})
\begin{align*}
X_i\ci X_j|\X_{\setminus{ij}}\Leftrightarrow J_{ij}=0.
\end{align*}
Then we define a GMRF as follows.

\begin{definition}
\label{def:gmrf}
\emph{(Definition 2.1, \citealp{RH:2005})} A random vector $\X$ is called a GMRF \emph{with respect to a graph} $\G=(\V,\E)$ with information matrix $J$ and potential vector $h$ if and only if its density has the form of (\ref{eq:cf}) and $J_{ij}\neq0\Leftrightarrow(i,j)\in\E$, for all $i$ and $j$.
\end{definition}
\noindent It is known that $X$ satisfies the Markov property on $\G$ (see Theorem 2.4 by \citealp{RH:2005})\footnote{For a GMRF, the pairwise Markov property, the local Markov property and the global Markov property are equivalent. This is proven by using Proposition 3.8 of \citet{Lau:1996}, in conjunction with the Hammersley-Clifford Theorem (Theorem 3.9, \citealp{Lau:1996}).}. Note that by Definition \ref{def:gmrf}, there is a one-to-one correspondence between the structure of $J$ and the structure of $\G$. Most information matrices $J$ for GMRFs are sparse, i.e.~there are $O(n)$ non-zero entries in $J$. The sparsity structure of $J$ facilitates the simulation of GMRFs through $J$ \citep{RH:2005}.

Another advantage of using the canonical form concerns the estimation of $\X$ given noisy observations $\mathbf{Y}$ \citep{Johnson:2002}. Indeed, assume that $\mathbf{Y}=C\X+\varepsilon$, where $C$ is an $n\times n$ real matrix and $\varepsilon\sim\N(0,M)$. Then the conditional distribution of $\X|\mathbf{Y}$ is given by 
\begin{align*}
f(x|y)&=\frac{f(y|x)f(x)}{f(y)}\\
&\propto f(y|x)f(x)\\
&\propto\exp\left[-\frac{1}{2}x^T\tilde{J}x+\tilde{h}^Tx\right],
\end{align*}
where $\tilde{J}=J+C^TM^{-1}C$ and $\tilde{h}=h+C^TM^{-1}y$. Thus given noisy observations, one only needs to update the information matrix and the potential vector to construct a graphical model for $f(x|y)$. For simplicity, we use $J$ and $h$ to denote the parameters after absorption of observations.

Given the canonical form, one needs both the covariance matrix $\Sigma=J^{-1}$ and the mean vector $\mu=J^{-1}h$ to obtain the marginal distributions of $\X$ or $\X|\mathbf{Y}$. Since knowing $J^{-1}$ is sufficient to recover $\mu$, we focus our efforts on calculating $J^{-1}$. Direct inversion of $J$ has complexity $O(n^3)$ and does not exploit the sparsity of $J$.\footnote{Algorithms that compute \emph{some} entries (either diagonal entries or certain off-diagonal entries) of a symmetric sparse matrix with complexity less than $O(n^3)$ do exist (see \citealp{EW:2013,LAKD:2008,LYMYE:2011,TS:2012}). The path-sum representation achieves this as well, computing the covariance of a pair of variables with complexity $O(n)$ whenever $\G$ is a tree: see \S\ref{CompCost} and \citet{Giscard2013}.} 
%Moreover, \emph{all} entries of $J^{-1}$ are needed to evaluate the mean vector $\mu$. 
In a simple situation where the graph $\G$ of a GMRF is a tree, belief propagation (BP) efficiently calculates the correct marginals \citep{Malioutov2006,Pearl:1988}. For  graphs with cycles (also called loopy graphs), the method of loopy belief propagation (LBP) can be used to efficiently \emph{approximate} the marginals. However, it was shown in \citet{WF:2001} that while LBP gives correct means for the marginals, the estimates of the covariance matrices it provides are generally incorrect.
%\pl{In the situation where the graph $\G$ of a GMRF is singly connected (i.e.~there is only one simple path between any pair of vertices),} as is the case of chains and trees, then \pl{belief propagation can be achieved exactly and} 

In this article, we present a novel approach to the calculation of the marginals of $\X$ or $\X|\mathbf{Y}$, which we term \emph{method of \,path-sums}. This approach is a generalization and completion of the walk-sum formulation developed in \citet{Malioutov2006}. The method of path-sums is based on results established by \citet{Giscard2012} concerning the algebraic structure of walk sets, which permit the systematic resummation of infinite families of walks in any walk-sum. These resummations transform a walk-sum into a branched continued fraction comprising only a \emph{finite} number of terms. Furthermore, these terms have an elementary  interpretation as simple paths and simple cycles on $\G$.  A simple path is a walk (i.e.~a trajectory on $\G$) whose vertices are all distinct. A simple cycle is a walk whose endpoints are identical and intermediate vertices (from the second to the penultimate) are all distinct and different from the endpoints. 
An important consequence of these observations is that if $J$ is positive definite\footnote{$J\succ0$ is equivalent to $J$ being non-singular for a GMRF.}, the path-sum formulation of $J^{-1}$ is convergent.
%, irrespectively of the walk-summability criterion that is so crucial to the walk-sum approach.

Consequently, one does not need the walk-summability of a model, which was devised by \cite{Malioutov2006} to guarantee the convergence of the infinite walk-sum representation of the covariance matrix $\Sigma=J^{-1}$. Let $R := I-J$ be the partial correlations matrix and $|R|$ be the entrywise absolute value of this matrix, i.e. $(|R|)_{ij} = |R_{ij}|$. The authors showed that a GMRF is walk-summable if and only if the spectral radius $\rho(|R|)$ is strictly less than 1, which implies $J\succ 0$. However, the converse does not hold: that is, there are positive definite information matrices which are not walk-summable. In contrast to the walk-summability criterion, our formulation only requires positive definiteness of $J$. Most importantly, the path-sum formulation gives the \emph{correct} marginals when the graph $\G$ is loopy, and is equivalent to BP when $\G$ is a tree.

The rest of this article is organized as follows. In \S\ref{sec:ps}, we introduce  the context and arguments underlying the path-sum formulation. We present the path-sum result in \S\ref{PSTHM} and show its validity for all positive definite matrices, irrespective of the walk-summability criterion. In the following section, \S\ref{ExistsSec}, we relate the path-sum representation of a covariance matrix to existing approaches. In \S\ref{BlockSec} we prove that the path-sum representation can be applied to arbitrary partitions of the information matrix $J$. We give an example demonstrating this claim. Finally in \S\ref{CompCost} we briefly discuss the computational cost of our approach as well as future prospects. The proof of the path-sum result is deferred to Appendix A.

\section{Path-sum representation}\label{sec:ps}
\subsection{Context}
The \textit{``most basic result of algebraic graph theory"} (as described in \citealp{Flajolet2009}) states that the powers of the adjacency matrix $A_\G$ of a graph $\G$ generate all the walks on this graph \citep{Biggs1993}. This result extends to weighted graphs if $A_\G$ is replaced by a weighted adjacency matrix $R$, with $R_{ij}$ the weight of the edge from vertex $j$ to vertex $i$ on $\G$. Then $(R^\ell)_{ij}$ is the sum of the weights of all the walks of length $\ell$ from $j$ to $i$ \citep{Flajolet2009}. The weight of a walk is simply the ordered product of the weight of the edges it traverses. (Note that the indices of a matrix are written right-to-left, but correspond to an edge written left-to-right. This is due to unfortunate conventions). We index the entries of $R$ by the labels of the vertices in $\G$. We write these labels with roman letters, Greek letters, or numbers, as convenient.
Now -- assuming that the Taylor series converges -- we have $(I-R)^{-1}=\sum_\ell R^\ell.$ It follows that $(I-R)^{-1}_{ij}$ can be interpreted as the sum of the weights of all the walks from $j$ to $i$. This directly implies the walk-sum interpretation advocated by \citet{Malioutov2006} for $J^{-1}=(I-R)^{-1}$, with $R=I-J$ and $\G$ the graphical model of the GMRF constructed by using Definition \ref{def:gmrf}. Further, it follows that the calculation of $J^{-1}$ is susceptible to a particular resummation technique from graph theory based on the structure of sets of walks, called the method of path-sums.

In its most general form, the method of path-sums stems from a fundamental algebraic property of the set of all walks on any weighted graph: namely, that any walk factorizes uniquely into products of prime walks, which are the simple paths and simple cycles of $\G$ (see \S\ref{sec:intro} for the definitions). The path-sum representation of the series of all walks on the graph $\G$ is thus the representation of this series that only involves the prime walks.\footnote{Since path-sums are the prime representations of walk series, they are the graph-theoretic analog of Euler product formulas for the Riemann zeta function and other totally multiplicative functions in number theory.} Since a finite graph sustains only finitely many primes, the walk series (which is typically infinite) thus has an exact representation involving only finitely many terms. An important consequence of this observation is that the path-sum expression of $J^{-1}$ is \emph{convergent} as long as $J\succ0$.

For a full exposition of the algebraic structure of walk sets at the origin of path-sums and its applications in linear algebra, the interested reader can refer to \citet{Giscard2012} and \citet{Giscard2013}. In \S\ref{PSTHM} we give the \emph{explicit} and universal path-sum formulation for $J^{-1}$. This expression takes the form of a branched continued fraction of finite depth and breadth.
% We begin by introducing the minimal notation required to present our results. 
%
%\subsection{Notation}
%\hspace{5mm}\textit{Graphs:} A graph $\G$ is a set of vertices connected by edges. These edges may be directed, in which case we say that $\G$ is a digraph, and we allow $\G$ to have self-loops (length 1 cycles). Following standard notations, we denote the vertices of $\G$ by Greek letters. We write $\G\backslash\{\alpha_1,\dotsc, \alpha_n\}$ for the graph obtained from $\G$ by deleting vertices $\alpha_1\dotsc \alpha_n$ as well as all the edges leading to and coming from these vertices.
%Consider a GMRF with information matrix $J$ and associated graphical model $\G$.
%Since $J$ and $\G$ have the same structure, we can associate to each edge of $\G$ a corresponding entry of $J$, called its weight. 

\subsection{Path-sum formulation of the covariance matrix}\label{PSTHM}
Let $\G\setminus\{\alpha,\beta\dotsc\}$ denote the subgraph of $\G$ obtained by deleting from $\G$ the vertices $\{\alpha,\beta\dotsc\}\subset \V$ and the edges incident to them. For simplicity, we write $J^{-1}_{\alpha\beta}$ for $(J^{-1})_{\alpha\beta}$.
The path-sum expression for $\Sigma=J^{-1}$ is presented in Theorem \ref{JPSresult} below. We defer its proof to Appendix A.
\begin{theorem}\label{JPSresult}
Let $J\succ0$ be an information matrix. Let $\Pi_{\G;\,\alpha\omega}$ and $\Gamma_{\G;\,\alpha}$ be the sets of simple paths from $\alpha$ to $\omega$ on $\G$ and the set of simple cycles from $\alpha$ to itself on $\G$, respectively. If $\G$ has finitely many vertices and edges, these two sets are finite.

Then each entry of the covariance matrix $\Sigma=J^{-1}$ admits an expression involving only weighted prime walks, called a path-sum representation. 
It is explicitly given by
\begin{align}
&\hspace{-1.5mm}J^{-1}_{\omega\alpha}=\sum_{p\in \Pi_{\G;\,\alpha\omega}}(-1)^{\ell(p)}\prod_{j=1}^{\ell(p)+1}\bigg\{\big(J_{\G\backslash\{\alpha,\nu_2,\dotsc,\nu_{j-1}\}}\big)^{-1}_{\nu_{j}\nu_{j}}\, J_{\nu_{j+1}\nu_{j}}\bigg\}\,J^{-1}_{\alpha\alpha}\,,\label{Jaw}\\
&\hspace{-1.5mm}J^{-1}_{\alpha\alpha}=\left(\sum_{\gamma\in\Gamma_{\G;\,\alpha\alpha}}\!(-1)^{\ell(\gamma)+1}\,J_{\mu_{1}\mu_{\ell(\gamma)}}\prod_{j=2}^{\ell(\gamma)}\bigg\{
\big(J_{\G\backslash\{\alpha,\mu_2,\dotsc,\mu_{j-1}\}}\big)^{-1}_{\mu_{j}\mu_{j}}\, J_{\mu_j\mu_{j-1}}\bigg\}\right)^{\!\!-1}\hspace{-3mm},\hspace{-1.5mm}\label{Jaa}
\end{align}
where the products are right-to-left (i.e.~$\prod_{i=1}^ma_i=a_m\cdots a_1$), $p=(\nu_1,\,\nu_2,\,\dotsc,\, \nu_{\ell(p)+1})$ is a simple path of length $\ell(p)$ with $\alpha\equiv \nu_1$ and $\omega\equiv \nu_{\ell(p)+1}$ for convenience; and $\gamma=(\mu_1,\,\mu_2,\,\dotsc,\, \mu_{\ell_\gamma},\,\mu_1)$ is a simple cycle of length $\ell(\gamma)$ from $\alpha\equiv \mu_1$ to itself. 
\end{theorem}

Note that $J^{-1}_{\alpha\alpha}$ is obtained recursively through Eq.~(\ref{Jaa}). Indeed it is expressed in terms of entries of inverses of submatrices of $J$, such as $(J_{\G\setminus\{\alpha,\mu_2,\dotsc, \mu_{j-1}\}})^{-1}_{\mu_j\mu_{j}}$, which is in turn obtained through Eq.~(\ref{Jaa}) but on the subgraph $\mathcal{G}\backslash{\{\alpha,\ldots,\mu_{j-1}}\}$ of $\mathcal{G}$. The recursion stops when vertex $\mu_j$ has no neighbor on this subgraph, in which case $(J_{\G\setminus\{\alpha,\mu_2,\dotsc, \mu_{j-1}\}})^{-1}_{\mu_j\mu_{j}}=1/J_{\mu_j\mu_j}$ (note that $J_{\mu_j\mu_j}\neq0$ since $J\succ0$). The entry $J^{-1}_{\alpha\alpha}$ is therefore expressed as a branched continued fraction which terminates at a finite depth, and $J^{-1}_{\alpha\omega}$ is a finite sum of such continued fractions. 

\begin{remark}\emph{
Theorem \ref{JPSresult} is valid even when $\G$ is loopy and/or $J$ is not walk-summable. In particular there is no restriction on the spectrum of $J$ as long as $J\succ0$. An example showing this is given below.}
\end{remark}

\begin{figure}[t!]
\vspace{-7mm}
\begin{center}
\includegraphics[width=.4\textwidth]{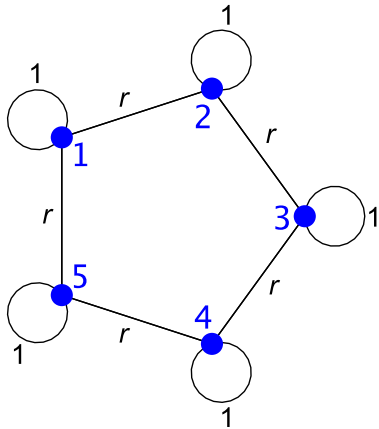}
\end{center}
\vspace{-5mm}
\caption{The circle graph on 5 vertices $\C_5$ \emph{associated with} $J$ in Example \ref{synthex}. The edge-weights are indicated next to the edges and the vertices are labeled 1 to 5.}\label{fig:C5}
\vspace{-1mm}
\end{figure}

\begin{example}[An illustrative example]\label{synthex}
\emph{To illustrate Theorem \ref{JPSresult}, we consider an example taken from \citet{Malioutov2006}, where $J$ has the structure of a circle graph on 5 vertices, denoted $\C_5$, see Fig.~\ref{fig:C5}. The information matrix is 
\begin{align*}
J=\begin{pmatrix}
1 & r & 0 & 0 & r\\
r & 1 & r & 0 & 0\\
0 & r & 1 & r & 0\\
0 & 0 & r & 1 & r\\
r & 0 & 0 & r & 1
\end{pmatrix}.
\end{align*}
Then $J$ is positive definite for $\frac{2}{1-\sqrt{5}}\leq r\leq\frac{2}{1+\sqrt{5}}$ (approximately $ -1.618\leq  r \leq 0.618$), and walk-summable if and only if $-1/2\leq r\leq 1/2$. This example thus provides a test case for both the walk-summable and non-walk-summable situations, depending on the value of $r$. Here we obtain $J^{-1}$ correctly for all values of $r$ such that $J$ is not singular.}
\end{example}

\vspace{-2mm}By the symmetry of $\C_5$, all the diagonal entries of $J^{-1}$ are identical, and we choose to calculate $J^{-1}_{11}$ without loss of generality. There are five simple cycles from vertex 1 to itself on $\C_5$: i) the self-loop $1\to1$; ii) the two backtracks $1\to 2\to 1$ and $1\to 5\to 1$; and iii) the two pentagons $1\to2\to3\to4\to5\to1$ and $1\to5\to4\to3\to2\to1$. By symmetry, the two backtracks and the two pentagons have the same weights.
Consequently Eq.~(\ref{Jaa}) gives
\begin{align*}
(J^{-1})_{11}=\bigg(&\underbrace{~~1\vphantom{\big)_{22}}~~}_{\text{Self-loop}\SelfLoopOne}~-~~\underbrace{~~2~ r^2\,\big(J_{\C_5\setminus\{1\}}\big)^{-1}_{22}\quad}_{\text{Backtracks }\BacktrackOne\BacktrackTwo}\\
&\vphantom{\Bigg{[}}\hspace{5mm}+~\underbrace{2\,r^5\big(J_{\C_5\setminus\{1,2,3,4\}}\big)^{-1}_{55}\big(J_{\C_5\setminus\{1,2,3\}}\big)^{-1}
_{44}\big(J_{\C_5\setminus\{1,2\}}\big)^{-1}_{33}\left(J_{\C_5\setminus\{1\}}\right)^{-1}_{22}}_{\text{Pentagons }\PentagonOne\PentagonTwo}~\bigg)^{-1}.
\end{align*}
The required entries $(J_{\C_5\setminus\{1\}})^{-1}_{22},\,\cdots,\,(J_{\C_5\setminus\{1,2,3,4\}})^{-1}_{55}$ remain to be calculated. To this end, we use again Eq.~(\ref{Jaa}) of Theorem \ref{JPSresult}. For example, consider calculating $(J_{\C_5\setminus\{1\}})^{-1}_{22}$. Since the graph associated with $J_{\C_5\setminus\{1\}}$ is $\C_5$ with vertex 1 removed, the only simple cycles from vertex 2 to itself on $\C_5\backslash\{1\}$ are the self-loop $2\to2$ and the backtrack $2\to 3\to2$. We thus find
\begin{equation*}
\left(J_{\C_5\setminus\{1\}}\right)^{-1}_{22}=\bigg(\underbrace{~~1\vphantom{\big)_{33}}~~}_{\text{Self-loop}\SelfLoopThree}-~\quad\underbrace{~r^2\big(J_{\C_5\setminus\{1,2\}}\big)^{-1}_{33}~}_{\text{Backtrack}\BacktrackThree}\bigg)^{-1},
\end{equation*}
Similarly we obtain $(J_{\C_5\setminus\{1,2\}})^{-1}_{33}=\big(1-r^2(J_{\C_5\setminus\{1,2,3\}})^{-1}_{44}\big)^{-1}$, $(J_{\C_5\setminus\{1,2,3\}})^{-1}_{44}=\big(1-r^2(J_{\C_5\setminus\{1,2,3,4\}})^{-1}_{55})^{-1}$ and finally $(J_{\C_5\setminus\{1,2,3,4\}})^{-1}_{55}=1$.
Combining these equations gives
\begin{align*}
(J_{\C_5\setminus\{1,2,3,4\}})^{-1}_{55}=1&\Rightarrow \left(J_{\C_5\setminus\{1,2,3\}}\right)^{-1}_{44}=\frac{1}{1-r^2}\\
&\hspace{8mm}\Rightarrow\left(J_{\C_5\setminus\{1,2\}}\right)^{-1}_{33}=\frac{1-r^2}{1-2r^2}\\
&\hspace{16mm}\Rightarrow\left(J_{\C_5\setminus\{1\}}\right)^{-1}_{22}=\frac{1-2r^2}{1-3r^2+r^4},\\
&\hspace{24mm}\Rightarrow(J^{-1})_{11}=\frac{1-3r^2+r^4}{1-5r^2+5r^4+2r^5}.
\end{align*}
As noted, the last equation gives the value of every diagonal entry of $J^{-1}$. We now consider the off-diagonal entries. We note again that, by the symmetry of $\mathcal{C}_5$, there are only two different entries, and we choose to calculate $(J^{-1})_{21}$ and $(J^{-1})_{31}$ without loss of generality. Following  Theorem \ref{JPSresult}, these two entries are given by a sum over simple paths from vertex 1 to vertex 2, and vertex 1 to vertex 3, respectively. In each case there are only two simple paths: for example from 1 to 2, we have $1\to2$ and $1\to5\to4\to3\to2$. Then Eq.~(\ref{Jaw}) yields
\begin{align*}
\big(J^{-1}\big)_{21}=&-\underbrace{~r\,\big(J_{\C_5\setminus\{1\}}\big)^{-1}_{22}\,\big(J^{-1}\big)_{11}~}_{\text{Simple path}\PathOnePentagon}\\
&+\underbrace{~r^4\,\big(J_{\C_5\setminus\{1,5,4,3\}}\big)^{-1}_{22}\,\big(J_{\C_5\setminus\{1,5,4\}}\big)^{-1}_{33}\,\big(J_{\C_5\setminus\{1,5\}}\big)^{-1}_{44}\,\big(J_{\C_5\setminus\{1\}}\big)^{-1}_{55}\,(J^{-1})_{11}~}_{\text{Simple path}\PathTwoPentagon}.
\end{align*}
By the symmetry of $\mathcal{C}_5$, we have $(J_{\C_5\setminus\{1,5,4,3\}})^{-1}_{22}=(J_{\C_5\setminus\{1,2,3,4\}})^{-1}_{55}$, $(J_{\C_5\setminus\{1,5,4\}})^{-1}_{33}=(J_{\C_5\setminus\{1,2,3\}})^{-1}_{44}$, $(J_{\C_5\setminus\{1,5\}})^{-1}_{44}=(J_{\C_5\setminus\{1,2\}})^{-1}_{33}$ and $(J_{\C_5\setminus\{1\}})^{-1}_{55}=(J_{\C_5\setminus\{1\}})^{-1}_{22}$. The previously obtained results then immediately give
\begin{align*}
&(J^{-1})_{21}=\frac{r^4+2r^3-r}{1-5r^2+5r^4+2r^5},\\
\shortintertext{and similarly}
&(J^{-1})_{31}=\frac{r^2-r^3-r^4}{1-5r^2+5r^4+2r^5}.
\end{align*}
Piecing these results together, we finally obtain
\begin{equation*}
J^{-1}=\frac{r^2}{2 r^3+3 r^2-r-1}\begin{pmatrix}\frac{(r-1) r-1}{r^2} & 1+r^{-1} & -1 & -1 & 1+r^{-1} \\
 1+r^{-1} & \frac{(r-1) r-1}{r^2} & 1+r^{-1} & -1 & -1 \\
 -1 & 1+r^{-1} & \frac{(r-1) r-1}{r^2} & 1+r^{-1} & -1 \\
 -1 & -1 & 1+r^{-1} & \frac{(r-1) r-1}{r^2} & 1+r^{-1} \\
 1+r^{-1} & -1 & -1 & 1+r^{-1} & \frac{(r-1) r-1}{r^2} 
 \end{pmatrix}.
\end{equation*}
We now verify that this expression is correct in both the walk-summable and non-walk-summable situations:
\\\\
\textbf{(a) Walk-summable and $J\succ 0$}: set $r=0.3$. Then to 5 decimal places we have $(J^{-1})_{11}=1.23975$, $(J^{-1})_{21}=-0.39959$ and $(J^{-1})_{31}=0.09221$, so that
\begin{align*}
J^{-1}=\begin{pmatrix}
1.23975 & -0.39959 & 0.09221 & 0.09221 & -0.39959\\
-0.39959 & 1.23975 & -0.39959 & 0.09221 & 0.09221\\
0.09221 & -0.39959 & 1.23975 & -0.39959 & 0.09221\\
0.09221 & 0.09221 & -0.39959 & 1.23975 & -0.39959\\
-0.39959 & 0.09221 & 0.09221 & -0.39959 & 1.23975
\end{pmatrix}.
\end{align*}
One can verify that this result obtained by using the path-sum formulation coincides with the one obtained by direct inversion of $J$.\\\\
\textbf{(b) Non walk-summable and $J\succ 0$}: set $r=0.6$. Then to 5 decimal places, we have $(J^{-1})_{11}=14.09091$, $(J^{-1})_{21}=-10.90909$ and $(J^{-1})_{31}=4.09091$, and so 
\begin{align*}
J^{-1}=\begin{pmatrix}
14.09091 & -10.90909 & 4.09091 & 4.09091 & -10.90909\\
-10.90909 & 14.09091 & -10.90909 & 4.09091 & 4.09091\\
4.09091 & -10.90909 & 14.09091 & -10.90909 & 4.09091\\
4.09091 & 4.09091 & -10.90909 & 14.09091 & -10.90909\\
-0.39959 & 4.09091 & 4.09091 & -0.39959 & 14.09091
\end{pmatrix}.
\end{align*}
This result is again easily verified through direct inversion of $J$.

\begin{remark}\emph{
First, one could have calculated everything with numerical values from the start, as opposed to evaluating the analytic expression of $J^{-1}$ for a specific value of $r$ as we did here. This gives the same results, as expected. Second, the analytical formula for $J^{-1}$ remains valid even when $J$ is not positive definite and only fails for those values of  $r$ such that $J$ is singular. The main \emph{mathematical} role of the condition $J\succ0$ in the path-sum representation is to guarantee that $J$ is not singular.\footnote{A path-sum result for singular matrices also exists, but it necessitates additional mathematical machinery. In this case the path-sum formulation yields a pseudo-inverse for the singular matrix, see \citet{Giscard2013}.}}
\end{remark} 

In the following section, we discuss the relations between Theorem \ref{JPSresult} and two existing approaches. 

\section{Relation to existing approaches}\label{ExistsSec}
\citet{Malioutov2006} provided a walk-sum derivation for Gaussian belief propagation on trees. \citet{Jones2005} presented an expression for the entries of $J^{-1}$ as a sum of simple paths. In this section, we show that these results are corollaries of Theorem \ref{JPSresult} arising as special cases. 

\subsection{Path-sums on trees}
The recursive structure of the path-sum representation of $J^{-1}$ is especially simple on trees, for which we recover the Gaussian belief propagation results of \citet{Malioutov2006}. 

Let $J\succ0$ be an information matrix associated with a tree model $\mathcal{T}$. 
For any vertex $\alpha$ of $\mathcal{T}$, let $\mathcal{N}(\alpha)$ be the set of neighbors of $\alpha$ on $\mathcal{T}$. Observe that since $\mathcal{T}$ is a tree, the only simple cycles from $\alpha$ to itself are the self-loop $\alpha\to\alpha$ with weight $J_{\alpha\alpha}$, and the backtracks to the neighbors of $\alpha$ on $\mathcal{T}$, e.g.~$\alpha\to\beta\to\alpha$, $\beta\in\mathcal{N}(\alpha)$. Then Eq.~(\ref{Jaa}) of Theorem \ref{JPSresult} gives
\begin{subequations}
\begin{equation}\label{stepMalioutov}
\Sigma_{\alpha\alpha}=J^{-1}_{\alpha\alpha}=\bigg(J_{\alpha\alpha}+\sum_{\beta\in\mathcal{N}(\alpha)}-J_{\alpha\beta}(J_{\mathcal{T}\setminus\{\alpha\}})^{-1}_{\beta\beta}J_{\beta\alpha}\bigg)^{-1}.
\end{equation}
The quantity $(J_{\mathcal{T}\setminus\{\alpha\}})^{-1}_{\beta\beta}$ satisfies a similar relation on the subtree $\mathcal{T}\backslash \{\alpha\}$, 
\begin{equation*}
(J_{\mathcal{T}\setminus\{\alpha\}})^{-1}_{\beta\beta}=\bigg(J_{\beta\beta}+\sum_{\delta\in\mathcal{N}(\beta)\backslash\alpha}-J_{\beta\delta}(J_{\mathcal{T}\setminus\{\alpha,\beta\}})^{-1}_{\delta\delta}J_{\delta\beta}\bigg)^{-1}.
\end{equation*}
Let $\C_{\mathcal{T}\backslash\{\beta\};\,\delta}$ be the connected component of $\mathcal{T}\backslash\{\beta\}$ that contains the vertex $\delta\in\mathcal{N}(\alpha)$. Similarly define $\C_{\mathcal{T}\backslash\{\alpha,\beta\};\,\delta}$ to be the connected component of $\mathcal{T}\backslash\{\alpha,\beta\}$ that contains the vertex $\delta$. Since $\mathcal{T}$ is a tree, these components are identical: $\C_{\mathcal{T}\backslash\{\beta\};\,\delta}=\C_{\mathcal{T}\backslash\{\alpha,\beta\};\delta}$, and therefore $(J_{\mathcal{T}\setminus\{\beta\}})^{-1}_{\delta\delta}=(J_{\mathcal{T}\setminus\{\alpha,\beta\}})^{-1}_{\delta\delta}$. Thus, we have
\begin{equation}\label{stepMalioutov2}
(J_{\mathcal{T}\setminus\{\alpha\}})^{-1}_{\beta\beta}=\bigg(J_{\beta\beta}+\sum_{\delta\in\mathcal{N}(\beta)\backslash\alpha}-J_{\beta\delta}(J_{\mathcal{T}\setminus\{\beta\}})^{-1}_{\delta\delta}J_{\delta\beta}\bigg)^{-1}.
\end{equation}
\end{subequations}
In order to show that Eqs.~(\ref{stepMalioutov}, \ref{stepMalioutov2}) are the Gaussian belief propagation results,  we introduce some notation from \citet{Malioutov2006}. 
%Let $r_{\beta\alpha}:=R_{\beta\alpha}=-J_{\beta\alpha}$ (see Eq.~(\ref{edgeweights})) and define 
%$\Delta J_{\beta\to\alpha} :=-J_{\alpha\beta}(J^{\alpha})^{-1}_{\beta\beta}J_{\beta\alpha}$, 
Let $\hat{J}_\alpha:=1/\Sigma_{\alpha\alpha}$ and $\hat{J}^{-1}_{\beta\backslash\alpha}:=(J_{\mathcal{T}\setminus\{\alpha\}})^{-1}_{\beta\beta}$.
%and $\gamma_{\beta\backslash\alpha}:=(J^{\alpha})^{-1}_{\beta\beta}$. 
With these notations, Eqs.~(\ref{stepMalioutov}, \ref{stepMalioutov2}) become
%J^{-1}_{\alpha\alpha}=J_{\alpha\alpha}+\sum_{\beta\in\mathcal{N}(\alpha)} \Delta J_{\beta\to\alpha}
\begin{align*}
&\hat{J}_{\alpha}=J_{\alpha\alpha}+\sum_{\beta\in\mathcal{N}(\alpha)}\Delta J_{\beta\to\alpha},\quad\text{and}~\quad\hat{J}_{\beta\backslash\alpha}=J_{\beta\beta}+\sum_{\delta\in\mathcal{N}(\beta)\backslash\alpha} \Delta J_{\delta\to\beta},
\end{align*}
with
\begin{equation*}
\Delta J_{\delta\to\beta}=-J_{\beta\delta}\hat{J}_{\delta\backslash\beta}^{-1}J_{\delta\beta}.
\end{equation*}
These are the Gaussian belief propagation equations (Eqs.~7, 8 and 9 in \citealp{Malioutov2006}), which immediately imply Propositions 16 and 17 in \citet{Malioutov2006} with the further definitions $r_{\beta\alpha}:=R_{\beta\alpha}=-J_{\beta\alpha}$, $\gamma_{\beta\backslash\alpha}:=(J_{\mathcal{T}\setminus\{\alpha\}})^{-1}_{\beta\beta}$. Note that $r_{\alpha\beta}=r_{\beta\alpha}$ since $J$ is symmetric. 
%
%In light of this derivation starting from Theorem \ref{JPSresult}, the path-sum representation of $J^{-1}$ appears as the generalization of the exact equations for Gaussian belief propagation on trees to all (loopy) graphs.

\subsection{An approach using determinants}
A determinant-based approach to the calculation of the covariance matrix $\Sigma=J^{-1}$ was demonstrated by \citet{Jones2005}. Here we show that this result follows from Eq.~(\ref{Jaw}) by using the adjugate formula for the matrix inverse. We then point out the fundamental limitation of determinant-based approaches. We overcome this limitation using the path-sum formulation in  \S\ref{BlockSec}.

\vspace{1mm}We recall the adjugate formula for matrix inverses:  
\begin{proposition}[Adjugate formula, \citealp{Strang2005}]\label{prop:adjugate}
Let $M\in\mathbb{C}^{n\times n}$ be a $n\times n$ non-singular complex matrix. Then $$\big(M^{-1}\big)_{ij}=(-1)^{i+j}\frac{\det M^{i}_{j}}{\det M},~\text{and in particular}~~\big(M^{-1}\big)_{ii}=\frac{\det M^{i}}{\det M},$$
where $M^i_j$ is the matrix $M$ with its $i$th column and $j$th row removed and $M^i=M_i^i$.
\end{proposition}
\noindent To obtain this result we start with Eq.~(\ref{Jaw}), which gives here
\begin{equation*}
J^{-1}_{\omega\alpha}=\sum_{p\in \Pi_{\G;\,\alpha\omega}}\!\!(-1)^{\ell(p)}(J_{\G\setminus\{\alpha,\,\nu_2,\dotsc,\nu_{\ell(p)}\}})^{-1}_{\omega\omega}J_{\omega\nu_{\ell(p)}}\dotsc J_{\nu_3\nu_2}(J_{\G\setminus\{\alpha\}})^{-1}_{\nu_2\nu_2}J_{\nu_2\alpha}(J^{-1})_{\alpha\alpha},
\end{equation*}
with $p=(\alpha,\nu_2,\dotsc \nu_{\ell(p)}\omega)$ is a simple path of length $\ell(p)$.
Since all the entries of $J$ commute, we can reorganize the weights in the above expression to get
\begin{align}
J^{-1}_{\omega\alpha}&=\sum_{p\in \Pi_{\G;\,\alpha\omega}}\!\!(-1)^{\ell(p)}\phi[p]\times (J_{\G\setminus\{\alpha,\,\nu_2,\dotsc,\nu_{\ell(p)}\}})^{-1}_{\omega\omega}\times\dotsc\times(J_{\G\setminus\{\alpha\}})^{-1}_{\nu_2\nu_2}\times(J^{-1})_{\alpha\alpha},\nonumber\\
\shortintertext{where $\phi[p]=J_{\omega\nu_{\ell(p)}}\dotsc J_{\nu_3\nu_2}J_{\nu_2\alpha}$ is the weight of the simple path $p$. \vphantom{\Big[}Using the adjugate formula in Proposition \ref{prop:adjugate}, we have}
J^{-1}_{\omega\alpha}&=\sum_{p\in \Pi_{\G;\,\alpha\omega}}\!\!(-1)^{\ell(p)}\phi[p]\times \frac{\det J^{\alpha,\,\nu_2,\dotsc,\,\nu_{\ell(p)},\,\omega}}{\det J^{\alpha,\,\nu_2,\dotsc,\,\nu_{\ell(p)}}}\times \dotsc \frac{\det J^{\alpha,\nu_2}}{\det J^\alpha}\times \frac{\det J^\alpha}{\det J},\nonumber\\
\shortintertext{since $J_{\G\setminus\{\alpha,\,\nu_2,\dotsc,\,\nu_{j}\}}=J^{\alpha,\,\nu_2,\dotsc,\,\nu_{j}}$. This simplifies to}
J^{-1}_{\omega\alpha}&=\sum_{p\in \Pi_{\G;\,\alpha\omega}}\!\!(-1)^{\ell(p)}\phi[p]\, 
\frac{\det J^{\alpha,\,\nu_2,\dotsc,\,\omega}}{\det J}.\label{JW2005}
\shortintertext{The last equation is given as Theorem 1 in \citet{Jones2005}.
Note that we only used Eq.~(\ref{Jaw}) in Theorem \ref{JPSresult} to obtain this result. \citet{Jones2005} did not give an expression equivalent to Eq.~(\ref{Jaa}) in Theorem \ref{JPSresult}, which provides an explicit formula for the necessary determinants in terms of simple cycles on the graph associated with $J$.}\nonumber
\end{align}

\vspace{-8mm}
The requirement that the entries of $J$ be commutative to obtain Eq.~(\ref{JW2005}) from Eq.~(\ref{Jaw}) may seem trivial, but in fact it excludes the very important case of \emph{block matrices}. In the next section, we demonstrate the use of Theorem \ref{JPSresult} in the case where $J$ is a block matrix.

\section{Path-sums for arbitrary partitions of $J$}\label{BlockSec}
\begin{figure}[t!]
\vspace{-5mm}
\begin{center}
\includegraphics[width=.9\textwidth]{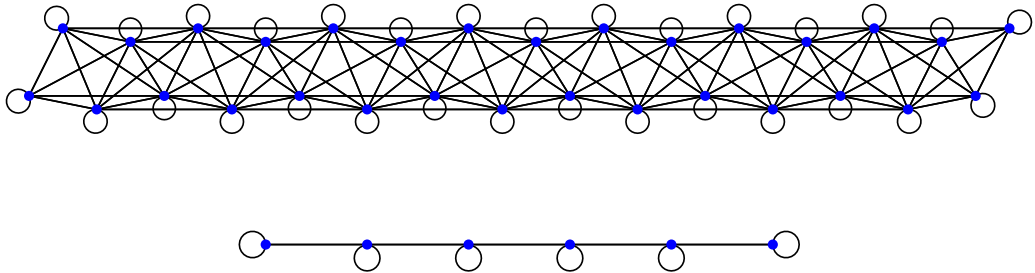}
\end{center}
\vspace{-2mm}
\caption{top: the graphical model associated with a 5-banded $30\times 30$ positive definite information matrix. The weights of the edges are the entries of $J$.
Bottom: the graphical model associated with a partition of $J$ into $5\times 5$ blocks. The edge weights are the blocks constituting $J$.
}\label{BlockBand}
\end{figure}

An information matrix $J$ may have a sparsity structure that is best exploited by partitioning $J$ into blocks. For example, consider an $n\times n$ information matrix $J$ which is banded: i.e.~$J_{ij}=0$ if and only if $|i-j|>b$ for some $b< n$. We say that $J$ is $b$-banded. Then $J$ is simply block-tridiagonal when partitioned $b\times b$ blocks. Consequently, the graph $\G$ associated with this partition of $J$ is simpler than the graph associated with the full matrix: see for example Fig.~\ref{BlockBand}. 
It is therefore desirable to develop a method capable of exploiting these simplifications.

A fundamental impediment to determinant-based approaches to $J^{-1}$ is that the notion of determinant itself does not extend to matrices with non-commutative blocks \citep{Silvester2000}.
%This observation stimulated research in non-commutative generalizations of the determinant in the 19th and 20th century, see e.g.~\citet{Gelfand2005}. 
%An information matrix may have a sparsity structure best exploited using blocks. 
%For example, consider the case where $J$ is a large $n\times n$, $n\gg 1$ banded matrix, i.e.~$J_{ij}=0$ if $|i-j|>b$ for some $b\leq n$. Then $J$ is also 
%%Easter egg simply block-tridiagonal. This is extremely important, especially in the case where $b\ll n$, since when seen as a block tridiagonal matrix, $J$ can be associated to a \emph{line graph} with $n/b$ vertices, much easier to deal with than the very complicated loopy graph obtained from the non-zero entries of $J$. 
In contrast, the path-sum formulation for $J^{-1}$ does not require the commutativity of the entries of $J$. For this reason, it continues to hold even when these entries do not commute, see \citep{Giscard2012} and \citep{Giscard2013}.  

Let $J$ be an $n\times n$ information matrix and $I_1,\dotsc,I_B$ for some $1\leq B\leq n$ be $B$ disjoint subsets of $\{1,\dotsc,n\}$. Then we write $J_{I_iI_j}$ for the minor (i.e.~block) of $J$ that corresponds to the rows indexed by $I_i$ and the columns indexed by $I_j$. We form a new graph $\G'=(\V',\E')$ associated with this partition of $J$, such that $J_{I_iI_j}\neq\mathbf{0}\iff (j,i)\in \E'$, where $\mathbf{0}$ is a zero matrix. The block $J_{I_iI_j}$ is now the weight associated with the edge $(j,i)$. With these conventions, Theorem \ref{JPSresult} extends to block matrices without modification.

\begin{remark}
\emph{For a GMRF $\X$, the partition of $J$ given above is equivalent to a partition of the set of random variables $\X$ into $B>1$ disjoint subsets $\X_1,\dotsc, \X_B$. Let $\X'=(\X_1,\dotsc,\X_B)$ be a GMRF with respect to a new graph $\G'=(\V',\E')$ with information matrix $J'$. Note that each $\X_i$ is now a random vector instead of a random variable. Following Definition \ref{def:gmrf}:
\begin{equation}\label{BlockGMRF}
\X_i\ci\X_j|\X'_{\setminus ij}\iff J'_{ij}=\mathbf{0}\iff (j,i)\notin\E',
\end{equation}
with $\X'_{\setminus ij}$ the set of variables with $\X_i$ and $\X_j$ removed from $\X'$. Note that $\X_i\ci\X_j|\X'_{\setminus ij}$ implies the global Markov property (GMP), which, for GMRFs, is equivalent to the pairwise Markov property (PMP) \citep{Lau:1996}. However for a general distribution, the PMP does not imply the GMP \citep{Lau:1996}. Hence in this case, the partition of $\X$ as per Eq.~(\ref{BlockGMRF}) may not correspond to the partition of $J$.}
\end{remark}
%We form a new graph $\tilde{\G}=(\tilde{\V},\tilde{\E})$ such that if for all $X\in \X_i$ and $X'\in \X_j$ we have $X\ci X'\,|\,\X_{X,\,X'}$ then $\{i,j\}\notin\tilde{\E}$.

Below we give a detailed synthetic example demonstrating the path-sum formulation of the covariance matrix $\Sigma=J^{-1}$ using a partition of $J$.

\begin{example}[Non walk-summable block matrix on a loopy graph]\label{blockex}
\emph{Consider the following positive definite information matrix $J$ of a thin membrane model 
\begin{align}\label{Jex}
&J=\\
&{\small \begin{pmatrix}
a+5 b & -b & 0 & -b & -b & 0 & -b & -b & 0 \\
 -b & a+5 b & -b & -b & 0 & -b & -b & 0 & -b \\
 0 & -b & a+5 b & 0 & -b & -b & 0 & -b & -b \\
 -b & -b & 0 & a+5 b & -b & 0 & -b & -b & 0 \\
 -b & 0 & -b & -b & a+5 b & -b & -b & 0 & -b \\
 0 & -b & -b & 0 & -b & a+5 b & 0 & -b & -b \\
 -b & -b & 0 & -b & -b & 0 & a+5 b & -b & 0 \\
 -b & 0 & -b & -b & 0 & -b & -b & a+5 b & -b \\
 0 & -b & -b & 0 & -b & -b & 0 & -b & a+5 b 
\end{pmatrix}},\nonumber
\end{align}
where $a,\,b>0$. Recall that walk-summability is equivalent to the condition $\rho(|R|)<1$, where $|R|$ is the entry-wise absolute value of $R:=I-J$ \citep{Malioutov2006}. Here we have $\rho(|R|)=\left| a+5 b-1\right| +\sqrt{19} \,b+b$, and thus -- depending on the values of $a$ and $b$ -- walk-summability does not always hold. For example, when $a=b=1$ we have $\rho(|R|)\simeq10.36$, $\rho(J)\simeq9.36$ and $\rho(R)\simeq8.36$ and $J$ is not walk-summable.}\end{example} 
\begin{figure}[t!]
\vspace{-5mm}
\begin{center}
\includegraphics[width=.8\textwidth]{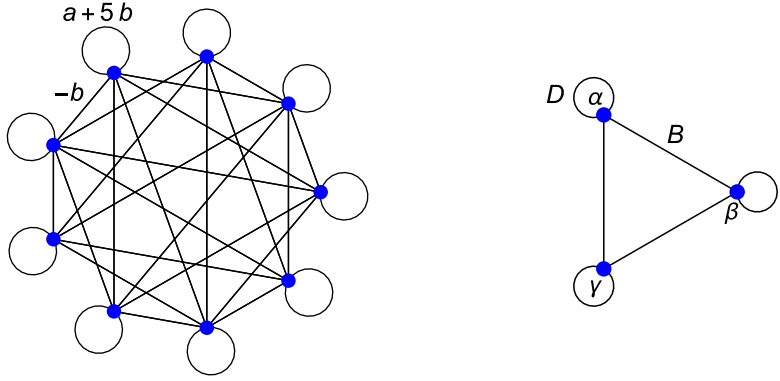}
\end{center}
\vspace{-9mm}
\caption{left: graphical model associated with the thin membrane model matrix $J$ of Eq.~(\ref{Jex}). Every edge has weight $-b$ and every self-loop has weight $a+5b$. Right: triangular graphical model $\mc{K}_3$ associated with the partition of $J$ given in Eq.~(\ref{JTriangle}). Edges and self-loops have matrix weights $E$ and $L$, respectively. The vertices are labeled $\alpha,~\beta$ and $\gamma$ for convenience.
}\label{exblock}
\vspace{-3mm}
\end{figure}

The graph $\G$ associated to $J$ is shown on the left of Fig.~\ref{exblock}. Instead of working on this complicated graph, we may partition the matrix into $3\times 3$ blocks as follows
\begin{align}\label{JTriangle}
J=\begin{pmatrix}L&E&E\\E&L&E\\E&E&L
\end{pmatrix},
\end{align}
with 
\begin{equation*}
E:=\begin{pmatrix} -b & -b & 0 \\
 -b & 0 & -b \\
 0 & -b & -b \end{pmatrix},\quad L:=\begin{pmatrix} a+5 b & -b & 0 \\
 -b & a+5 b & -b \\
 0 & -b & a+5 b\end{pmatrix}.
\end{equation*}
The graph $\G'$ associated with this partition of $J$ is a triangle, denoted $\mathcal{K}_3$ and shown on the right of Fig.~\ref{exblock}. Accordingly, the edge-weights are now the matrices $E$ and $L$ for edges and self-loops, respectively. Note that these weights do not commute: $[E,L]=EL-LE\neq 0$.

The only simple cycles from any vertex (say, $\alpha$) to itself on this graph are the self-loop $\alpha\to\alpha$, the two backtracks $\alpha\to\beta\to\alpha$ and $\alpha\to\gamma\to\alpha$ and the two triangles $\alpha\to\beta\to\gamma\to\alpha$ and $\alpha\to\gamma\to\beta\to\alpha$. Theorem \ref{JPSresult} therefore yields
\begin{subequations}
\begin{align}\label{Jaatriangle}
J^{-1}_{\alpha\alpha}=\bigg(&~\underbrace{~~L\vphantom{)_{\alpha}}~~}_{\text{Self loop}\LoopTri}~~-~~\underbrace{~E\,.\,(J_{\mathcal{K}_3\setminus\{\alpha\}})^{-1}_{\beta\beta}\,.\,E~}_{\text{Backtrack}\BackTriOne}~~-~~\underbrace{~E\,.\,(J_{\mathcal{K}_3\setminus\{\alpha\}})^{-1}_{\gamma\gamma}\,.\,E~}_{\text{Backtrack}\BackTriTwo}\\
&\hspace{-12mm}+~\underbrace{~E\,.\,(J_{\mathcal{K}_3\setminus\{\alpha,\beta\}})^{-1}_{\gamma\gamma}\,.\,E\,.\,(J_{\mathcal{K}_3\setminus\{\alpha\}})^{-1}_{\beta\beta}\,.\,E~}_{\text{Triangle}\TriOne}+\underbrace{~E\,.\,(J_{\mathcal{K}_3\setminus\{\alpha,\gamma\}})^{-1}_{\beta\beta}\,.\,E\,.\,(J_{\mathcal{K}_3\setminus\{\alpha\}})^{-1}_{\gamma\gamma}\,.\,E~}_{\text{Triangle}\TriTwo}\,\bigg)^{\!-1}\!\!,\nonumber
\end{align}
where, by Eq.~(\ref{Jaa}),
\begin{align}\label{Jawtriangle}
\hspace{-4mm}(J_{\mathcal{K}_3\setminus\{\alpha\}})^{-1}_{\beta\beta}=&\bigg(\hspace{-3.5mm}\underbrace{~~L\vphantom{)_{\gamma}}~~}_{\text{Self loop}\LoopTriNoAlpha}\hspace{-1.5mm}-\underbrace{~~E\,.\,(J_{\mathcal{K}_3\setminus\{\alpha,\beta\}})^{-1}_{\gamma\gamma}\,.\,E~~}_{\text{Backtrack}\BackTriNoAlpha}\bigg)^{-1}\!\text{ and }~
(J_{\mathcal{K}_3\setminus\{\alpha,\beta\}})^{-1}_{\gamma\gamma}=\hspace{-2mm}\underbrace{~L^{-1}\vphantom{)_{\alpha}}~}_{\text{Self loop}\LoopLoop}\hspace{-3mm}.
\end{align}
\end{subequations}
Since $(J_{\mathcal{K}_3\setminus\{\alpha\}})^{-1}_{\gamma\gamma}=(J_{\mathcal{K}_3\setminus\{\alpha\}})^{-1}_{\beta\beta}$ and $(J_{\mathcal{K}_3\setminus\{\alpha,\gamma\}})^{-1}_{\beta\beta}=(J_{\mathcal{K}_3\setminus\{\alpha,\beta\}})^{-1}_{\gamma\gamma}$, Eqs.~(\ref{Jaatriangle}, \ref{Jawtriangle}) completely determine $J^{-1}_{\alpha\alpha}$. Thus, we have
\begin{align*}
&J^{-1}_{\alpha\alpha}=\frac{b}{a^2+8 a b-3 b^2}\times\\
&\hspace{5mm}
\begin{pmatrix}
 \frac{a}{b}+\frac{a}{a+3 b}+\frac{5 a}{3 (a+6 b)}+\frac{1}{3} & 1 & \frac{3 b}{a+3 b} \\
 1 & \frac{a}{b}+\frac{12 a}{5 (a+5 b)}+\frac{3}{5} & 1 \\
 \frac{3 b}{a+3 b} & 1 & \frac{a}{b}+\frac{a}{a+3 b}+\frac{5 a}{3 (a+6 b)}+\frac{1}{3}
\end{pmatrix}.
\end{align*}
By symmetry of $J$, this gives also $J^{-1}_{\beta\beta}$ and $J^{-1}_{\gamma\gamma}$. The off-diagonal blocks, which are all identical because of the symmetry of the matrix, now follow easily. For example, we have
\begin{align*}
(J^{-1})_{\beta\alpha}&=\underbrace{\,(J_{\mathcal{K}_3\setminus\{\alpha\}})^{-1}_{\beta\beta}\,.\,E\,.(J^{-1})_{\alpha\alpha}\,}_{\text{Path}\PathOneTri}~-~\underbrace{E\,.\,(J_{\mathcal{K}_3\setminus\{\alpha,\gamma\}})^{-1}_{\beta\beta}\,.\,E\,.\,(J_{\mathcal{K}_3\setminus\{\alpha\}})^{-1}_{\gamma\gamma}\,.\,E\,.(J^{-1})_{\alpha\alpha}}_{\text{Path}\PathTwoTri},\nonumber\\
&\nonumber\\
&=\frac{b}{a^2+8 a b-3 b^2}\begin{pmatrix}  \frac{a}{a+3 b}-\frac{5a}{6 (a+6 b)}+\frac{5}{6} & 1 & \frac{3 b}{a+3 b} \\
 1 & \frac{6 b}{a+5 b} & 1 \\
 \frac{3 b}{a+3 b} & 1 & \frac{a}{a+3 b}-\frac{5a}{6 (a+6 b)}+\frac{5}{6} 
 \end{pmatrix}.
\end{align*} 
\noindent As long as $a,b>0$ (which implies $J\succ0$), these results hold regardless of the values of $a$ and $b$. For example setting $a=b=1$, we obtain the correct covariance matrix even though $J$ is not walk-summable. 

\section{Conclusion}\label{CompCost}
We have presented an approach for exact inference of marginals in Gaussian graphical models, called the path-sum formulation. This approach represents the covariance matrix $\Sigma$ of a Gaussian Markov Random Field (GMRF) as a continued fraction of finite depth and breadth, and which only involves the simple paths and simple cycles on the graphical model of the GMRF. For GMRFs, we have shown that the path-sum formulation only requires $J\succ0$. It also extends to arbitrary partitions of $J$, thereby providing the flexibility necessary to exploit the sparsity structure of the information matrix $J=\Sigma^{-1}$.

From an algorithmic point of view, the computational cost associated with the calculation of any entry of a $n\times n$ covariance matrix using a path-sum is known to be $O(n)$ on trees, see \cite{Giscard2013}. The $O(n)$ complexity is achieved for the most general situation; if $J$ or $\G$ has symmetries, the complexity would be lower. Evaluating the complexity of the path-sum approach on a graph with an arbitrary topology remains an open problem. Since the method of path-sums arises from exact resummations on walk-sums, every problem that has a walk-sum interpretation is susceptible to these resummations and therefore admits a path-sum expression. Consequently, every algorithm that has a walk-sum interpretation (see e.g.~\citealp{Chandrasekaran2008}) necessarily has a path-sum formulation.

The interpretation of the entries of the covariance matrix as walk-sums and the existence of the path-sum formulation opens the door to many more walk-based methods, which are intermediary between these two approaches. The central idea is to exploit recent results regarding the algebraic structure of the set of walks on graphs with arbitrary topology \citep{Thwaite2014} to identify certain infinite geometric series of terms appearing in a walk-sum. These geometric series can then be exactly resummed, thereby reducing the sum of all walk weights to a sum over the weights of a certain (yet infinite) subset of `irreducible' walks. Each term in this sum is `dressed' so as to exactly include the contributions of the infinite families of resummed terms. The exact form of both the dressing and the irreducible terms remaining in the sum depend on the structure of the resummed terms: choosing a different family of terms produces a different series. 

Viewed in this context, the walk-sum and path-sum formulations of a given problem can be seen to be the extrema of a hierarchy of possible resummations: a walk-sum corresponds to the case where \emph{no} terms are exactly resummed, such that an explicit summation over all walks remains to be carried out, while the path-sum expression corresponds to the case where \emph{all possible} geometric series have been resummed, leaving behind only a (finite) sum over simple paths. In between these extremes lie many intermediate formulations, whose mathematical properties such as complexity and convergence are expected to interpolate between those of walk-sums and path-sums.

\acks{P.-L.~Giscard and D.~Jaksch acknowledge funding from EPSRC Grant EP/K038311/1. D.~Jaksch also received funding from the ERC under the European Unions Seventh Framework Programme (FP7/2007-2013)/ERC Grant Agreement no.~319286 Q-MAC. Z.~Choo was funded by a departmental studentship of the Department of Statistics, University of Oxford. S.~J.~Thwaite acknowledges funding from the Alexander von Humboldt Foundation.}

\appendix
\section*{Appendix A.}
In this appendix we prove Theorem \ref{JPSresult}.\\

\begin{proof}
The proof of the theorem is organized in two steps. First, assuming walk-summability, we obtain the path-sum formulation for $J^{-1}$. Second, we show that path-sum expression thus obtained is the unique analytic continuation of the sum of all walks, and continues to exist in the absence of walk-summability. Consequently, the path-sum remains a valid representation of $J^{-1}$ when walk-summability fails.
%For the sake of completeness, we include a brief presentation of the (well established) theory underlying analytic continuation.  

\emph{Step 1:} the path-sum expression for $J^{-1}$ is a special case of the more general result concerning the path-sum formulation of the matrix inverse function presented and proved in \citet{Giscard2013}. Below we briefly recount how we obtain result for the specific case of $J^{-1}$.
Let $R:=I-J$ and $|R|$ be the entry-wise absolute value of $R$, i.e.~$|R|_{ij}=|R_{ij}|$. Let $\rho(|R|)$ be the spectral radius of $|R|$, i.e.~the largest eigenvalue of $|R|$. As established by \citet{Malioutov2006}, $\rho(|R|)<1$ is equivalent to to walk-summability, and is enough to guarantee absolute convergence of the series $\sum_{n\geq 0} R^n=(I-R)^{-1}=J^{-1}$. This power series can be seen as a sum of walk weights on the graph $\G$ associated with $R$ \citep{Flajolet2009}, that is
\begin{equation}\label{SumWalksR}
\Big(\sum_{n\geq 0}R^n\Big)_{\omega\alpha}=\sum_{w\in \W_{\G;\,\alpha\omega}}\phi[w],
\end{equation}
with $\phi[w]$ the weight of the walk $w$, which is defined to be the product of the weights of the edges traversed by $w$, where the weight of an edge from $\alpha$ to $\beta$ is given by
\begin{equation}\label{edgeweights}
\phi[\alpha\beta]:=R_{\beta\alpha}=\begin{cases}-J_{\beta\alpha}&\text{if }\alpha\neq \beta,\\
1-J_{\alpha\alpha}&\text{if }\alpha=\beta.
\end{cases}
\end{equation}
Note that since we assumed walk-summability, the right-hand side of Eq.~(\ref{SumWalksR}) exists.
We then use the result by \citet{Giscard2012}, which reduces a series of weighted walks, such as the one of Eq.~(\ref{SumWalksR}), to a sum of weighted simple paths and simple cycles. This result is reproduced here for the sake of completeness:
\begin{theorem}[Path-sum, \citealp{Giscard2012}]\label{PSformalresult}
Let $\G$ be a graph and $\phi[.]$ be the weight function of Eq.~(\ref{edgeweights}). Suppose that the walk-sum $\sum_{w\in \W_{\mathcal{G};\alpha\omega}}\phi[w]$ exists. Then this sum is given by the weighted path-sum
\begin{equation*}
\sum_{w\in \W_{\mathcal{G};\alpha\omega}}\!\!\!\!\phi[w]=\!\sum_{p\in \Pi_{\G;\,\alpha\omega}}\prod_{j=1}^{\ell(p)+1}\bigg\{\phi_{\G\backslash\{\alpha,\,\nu_2,\dotsc,\,\nu_{j-1}\};\,\nu_{j}}\, \phi[\nu_{j+1}\nu_{j}]\bigg\}\,\phi_{\G;\,\alpha},
\end{equation*}
where $\phi_{\G;\,\alpha}:=\sum_{w\in \W_{\mathcal{G};\alpha\alpha}}
\phi[w]$ is the weighted sum of all walks from $\alpha$ to itself on $\G$ and is explicitly given by
\begin{equation*}
\phi_{\G;\,\alpha}=\left(1-\!\!\!\sum_{\gamma\in\Gamma_{\G;\,\alpha\alpha}}\!\!\phi[\mu_{1}\mu_{\ell(\gamma)}]\prod_{j=2}^{\ell(\gamma)}\bigg\{
\phi_{\G\backslash\{\alpha,\,\mu_2,\dotsc,\,\mu_{j-1}\};\,\mu_{j}}\, \phi[\mu_j\mu_{j-1}]\bigg\}\right)^{-1}\!\!\!\!\!,
\end{equation*}
and similarly for all $\phi_{\G\backslash\{\alpha,\,\nu_2,\dotsc,\,\nu_{j-1}\};\,\nu_{j}}$ and $\phi_{\G\backslash\{\alpha,\,\mu_2,\dotsc,\,\mu_{j-1}\};\,\mu_{j}}$. In these expressions, the products are right-to-left, $p=(\nu_1,\,\nu_2,\,\dotsc,\, \nu_{\ell(p)+1})$ is a simple path of length $\ell(p)$ with $\alpha\equiv \nu_1$ and $\omega\equiv \nu_{\ell(p)+1}$ for convenience; and $\gamma=(\mu_1,\,\mu_2,\,\dotsc,\, \mu_{\ell_\gamma},\,\mu_1)$ is a simple cycle of length $\ell(\gamma)$ from $\alpha\equiv \mu_1$ to itself.  
\end{theorem}
\noindent By using the edge weights of Eq.~(\ref{edgeweights}), we obtain $\phi_{\G;\,\alpha}=J_{\alpha\alpha}^{-1}$ and similarly $\phi_{\G\backslash\{\alpha,\,\nu_2,\dotsc,\,\nu_{j-1}\};\,\nu_{j}}=(J_{\G\setminus\{\alpha,\,\nu_2,\dotsc,\,\nu_{j-1}\}})^{-1}_{\nu_j\nu_j}$ and Theorem \ref{PSformalresult} yields Eqs.~(\ref{Jaw}, \ref{Jaa}).  
%%\textcolor{red}{
%Remark that a self-loop is a length one simple cycle and is therefore present in the sum over simple cycles on the right-hand side of Eq.~(\ref{Saa}). Since $J\succ0$ all its diagonal entries are non-zero and are precisely the weights of these self-loops. According to Eq.~(\ref{Saa}), a self-loop contributes a weight of $\phi(\alpha\alpha)=R_{\alpha\alpha}=1-J_{\alpha\alpha}$ so that .}
Next we prove that Eqs.~(\ref{Jaw}, \ref{Jaa}) yield $J^{-1}$ for any matrix $J\succ0$, even when walk-summability does not hold and the walk-sum of Eq.~(\ref{SumWalksR}) does not converge. Central to our proof is the (well-established) theory of analytic continuation \citep{Priestley:2003}. 

\emph{Step 2:} we consider the three following functions of a complex variable $z$ into $\mathbb{C}^{n\times n}$, $g_1(z):=(I-zR)^{-1}$, $g_2(z):=\sum_{n\geq 0} z^n R^n$ and, since $g_2(z)$ is a walk-sum, it has a path-sum expression which we denote $g_3(z)$. (The path-sum expression $g_3(z)$ is obtained from Theorem \ref{JPSresult} on using the edge-weights $zR_{\beta\alpha}$ for an edge from $\alpha$ to $\beta$). 

The first function, $g_1(z)$, is analytic on $z\in \mathbb{C}\backslash \text{Sp}^{-1}(R)$ with $\text{Sp}^{-1}(R)$ the inverse of the spectrum of $R$. The second function, $g_2(z)$, is analytic on the disk $D_R$ of the complex plane where $z< 1/\rho(|R|)$. The third function, $g_3(z)$, comprises only a finite number of terms (since there are finitely many simple paths and simple cycles on a finite graph). Consequently $g_3(z)$ exists for $z\in\mathbb{C}\backslash \text{Sp}^{-1}(R)$. Evidently $g_1$ and $g_2$ agree on $D_R$ and by construction $g_2$ and $g_3$ agree on $D_R$ as well. It follows that $g_1$ and $g_3$ constitute two \emph{direct analytic continuations} of $g_2$ outside of $D_R$ \citep{Priestley:2003} (also called \emph{extensions} of $g_2$, \citealp{Kreyszig:1978}). By the Uniqueness Theorem (Theorem 15.9, \citealp{Priestley:2003}) only one such analytic continuation exists and $g_1(z)=g_3(z)$ on the domain $\mathbb{C}\backslash \text{Sp}^{-1}(R)$.

We conclude the proof by showing that the point $z=1$, for which $g_1(z)$ and $g_3(z)$ yield $J^{-1}$, is also in this domain. Indeed, the covariance matrix $\Sigma$ of a Gaussian distribution must be positive definite, implying that it is non-singular (Corollary 7.17 in \citealp{HJ:2013}). Consequently, $J$ is not singular: so $J^{-1}$ exists exists and $1$ is not an eigenvalue of $R=I-J$, i.e.~$1\in \mathbb{C}\backslash \text{Sp}^{-1}(R)$. Then $g_1$ and $g_3$ exist at $z=1$, in particular $g_3(1)$ is a valid representation of $J^{-1}$, even though $z=1$ may not be in $D_R$ (i.e.~$J$ is not walk-summable). This completes the proof of Theorem \ref{JPSresult}. 
\end{proof}

\vskip 0.2in
\bibliography{Stat_Inference}

\end{document}